\documentclass[english,final,leqno,onefignum,onetabnum]{siamltexmm}
\usepackage[T1]{fontenc}
\usepackage[latin9]{inputenc}
\usepackage{array}
\usepackage{booktabs}
\usepackage{hyperref}
\usepackage{textcomp}
\usepackage{multirow}
\usepackage{amsmath}
\usepackage{algorithm2e}
\usepackage{amssymb}
\usepackage{graphicx}
\usepackage{setspace}
\usepackage{esint}
\usepackage{babel}

\title{An Efficient Non-Intrusive Uncertainty Propagation Method for Stochastic Multi-Physics Models} 

\author{A.~Mittal\footnotemark[1] \footnotemark[3]
\and G. Iaccarino\footnotemark[2]}

\begin{document}
\maketitle
\newcommand{\slugmaster}{%
\slugger{juq}{xxxx}{xx}{x}{x--x}}

\renewcommand{\thefootnote}{\fnsymbol{footnote}}
\footnotetext[1]{Institute for Computational and Mathematical Engineering, Stanford University, Stanford, CA 94305.}
\footnotetext[2]{Mechanical Engineering, Stanford University, Stanford, CA 94305.}
\footnotetext[3]{Corresponding author (Email: \email{mittal0@stanford.edu})}
\renewcommand{\thefootnote}{\arabic{footnote}}

\begin{abstract}
Multi-physics models governed by coupled partial differential equation (PDE) systems, are naturally suited for partitioned, or modular numerical solution strategies. Although widely used in tackling deterministic coupled models, several challenges arise in extending the benefits of modularization to uncertainty propagation. On one hand, Monte-Carlo (MC) based methods are prohibitively expensive as the cost of each deterministic PDE solve is usually quite large, while on the other hand, even if each module contains a moderate number of uncertain parameters, implementing spectral methods on the combined high-dimensional parameter space can be prohibitively expensive. In this work, we present a reduced non-intrusive spectral projection (NISP) based uncertainty propagation method  which separates and modularizes the uncertainty propagation task in each subproblem using block Gauss-Seidel (BGS) techniques. 
The overall computational costs in the proposed method are also mitigated by constructing reduced approximations of the input data entering each module. These reduced approximations and the corresponding quadrature rules are constructed via simple linear algebra transformations. We describe these components of the proposed algorithm assuming a generalized polynomial chaos (gPC) model of the stochastic solutions. We demonstrate our proposed method and its computational gains over the standard NISP method using numerical examples.
\end{abstract}



\pagestyle{myheadings}
\thispagestyle{plain}
\markboth{A.~MITTAL AND G. IACCARINO}{NON-INTRUSIVE UNCERTAINTY PROPAGATION}

\section{Introduction}

With the aim of validating predictions from computer simulations against
real-world experiments, it becomes necessary to include uncertainties (descrbied as random quantities) within the
associated physical model. Uncertainties of this kind, known as aleatoric uncertainties, mainly arise
because mathematical models are often idealized approximations
of their target scenarios (in terms of geometry, material properties and boundary constraints
), and due to limited knowledge
in defining parameters of the physical systems being investigated. Uncertainty Quantification
(UQ) tools, have therefore, become key requirements for estimating
the credibility and confidence of predictions from computer simulations.
In particular, addressing the computationally demanding task of
uncertainty propagation, has gained tremendous prominence amongst simulation engineers, supported by the rapid growth in computing power and advances in numerical analysis.

To this end, Monte-Carlo (MC) based methods [\hyperlink{ref1}{1}, \hyperlink{ref2}{2}], wherein
the corresponding ouput samples for a set of random input samples
can be collected via repeated simulation runs, provide the most straightforward approach. 
However, when dealing with costly
PDE solvers, computing an ensemble of output realizations large enough
to evaluate accurate statistics can be prohibitively expensive. Alternatively,
cheaper surrogate models, can be trained on a much smaller ensemble
of solution samples, using various sparse regression techniques [\hyperlink{ref3}{3}, 
\hyperlink{ref4}{4}], and subsequently sampled exhaustively to approximate solution
distributions and/or statistics. However, in the context of multi-physics
problems, or coupled PDE systems in general, these methods ignore the
rich and potentially exploitable structures within the model.

Several coupled PDE systems, for example, fluid-structure [\hyperlink{ref5}{5}] and reactive-transport
models [\hyperlink{ref6}{6}], are natural targets for employing modular (partitioned)
solution methods wherein, solvers can be rapidly
developed by reusing legacy solvers for each
constituent subproblem. The disadvantage however, is that the convergence
rate of modular methods, can be much slower when compared to monolithic
methods, for example, the quadratically convergent Newton's method
[\hyperlink{ref7}{7}]. Therefore, in the context of coupled PDE solvers, modularization typically results in a trade-off
between the developmental and computational costs, of which the former
usually dominates. However, in extending the standard practices of modularization towards uncertainty propagation, several challenges yet remain.

In this work, we propose an efficient non-intrusive method for generalized
polynomial chaos (gPC) [\hyperlink{ref8}{8}] based uncertainty propagation in coupled
PDE solvers, which mitigates the curse of dimensionality associated
with spectral uncertainty propagation methods. This is achieved by constructing reduced dimensional (and order) approximations of the input data are before they enter their respective solver modules. Correspondingly,
reduced approximations of the respective output data are obtained
using optimal quadrature rules, which consequently reduce the required
number of repeated module runs significantly. Simple linear algebraic tools
are used in constructing the reduced approximations and quadrature
rules. Moreover, the approximation errors can be controlled separately
in each module. While several recent works have demonstrated the possibility
of exploiting relatively trivial coupling structures, for example, unidirectional
coupling [\hyperlink{ref9}{9}], and linear coupling [\hyperlink{ref10}{10}, \hyperlink{ref11}{11}, \hyperlink{ref12}{12}], to reduce the costs
of uncertainty propagation, we consider a more general setup with bidirectional and nonlinear coupling structures in this work. 

Our proposed method is a generalization of 
the approach recently proposed by Constantine et. al [\hyperlink{ref13}{13}], which
applies well to network (weakly) coupled multi-physics systems. Moreover, several
components of our proposed algorithm have also been motivated by the
recent works of Arnst et. al. [\hyperlink{ref}{14}, \hyperlink{ref15}{15}], wherein reduced
chaos expansions [\hyperlink{ref16}{16}] are used in approximating the input data. 
However, in their approach, the dimension reduction procedure can only be implemented unidirectionally, and therefore, the computational gains
achieved would be limited to those respective modules only. We address this limitation in our proposed method and facilitate the reduction strategy to be implemented across all the modules. Furthermore, 
the respective tolerances, which control the error in the reduced approximations, can be prescribed individually for each module.

The remainder of this article is organized as follows. In \hyperlink{sec2}{\S2}, we
provide an overview of the preliminary definitions and concepts, along with a description of the standard non-intrusive spectral projection (NISP) method for solving stochastic algebraic systems. 
In \hyperlink{sec3}{\S3}, we describe the proposed
reduced NISP method and its constituent
dimension and order reduction steps. Moreover, we prove the existence of asymptotic upper bounds
on the approximation errors incurred by our proposed method. In \hyperlink{sec4}{\S4}, we
report and compare the performance and accuracy of the standard and reduced NISP method implementations on
two numerical examples.

\section{Preliminary definitions and concepts}
\hypertarget{sec2}{}
Without loss of generality, we focus our analysis on a two component, steady-state
coupled PDE system. First, we consider the deterministic case.

\subsection{Deterministic model}

A spatial discretization of the system yields
the coupled algebraic system:
\hypertarget{eq21}{}
\begin{align}
 & \boldsymbol{f}_{1}\left(\boldsymbol{u}_{1};\boldsymbol{u}_{2}\right)=\boldsymbol{0},\ \boldsymbol{f}_{1},\boldsymbol{u}_{1}\in\mathbb{R}^{n_{1}},\nonumber \\
 & \boldsymbol{f}_{2}\left(\boldsymbol{u}_{2};\boldsymbol{u}_{1}\right)=\boldsymbol{0},\ \boldsymbol{f}_{2},\boldsymbol{u}_{2}\in\mathbb{R}^{n_{2}},
\end{align}
where $\boldsymbol{f}_{1}$, $\boldsymbol{f}_{2}$ denote the discrete component
residuals and $\boldsymbol{u}_{1}$, $\boldsymbol{u}_{2}$ denote
the respective finite-dimensional discretizations of the solution
fields intrinsic in each component PDE system.

\subsection{Deterministic numerical solvers}

There are various approaches to solve \hyperlink{eq21}{Eq. 2.1} numerically. Monolithic
solution methods, for example, Newton's method, would require repeatedly solving
the fully-coupled linear system
\hypertarget{eq22}{}
\begin{equation}
\left[\begin{array}{cc}
{\displaystyle \frac{\partial\boldsymbol{f}_{1}}{\partial\boldsymbol{u}_{1}}}\left(\boldsymbol{u}_{1}^{\ell};\boldsymbol{u}_{2}^{\ell}\right) & {\displaystyle \frac{\partial\boldsymbol{f}_{1}}{\partial\boldsymbol{u}_{2}}}\left(\boldsymbol{u}_{1}^{\ell};\boldsymbol{u}_{2}^{\ell}\right)\\
\\
{\displaystyle \frac{\partial\boldsymbol{f}_{2}}{\partial\boldsymbol{u}_{1}}}\left(\boldsymbol{u}_{2}^{\ell};\boldsymbol{u}_{1}^{\ell}\right) & {\displaystyle \frac{\partial\boldsymbol{f}_{2}}{\partial\boldsymbol{u}_{2}}}\left(\boldsymbol{u}_{2}^{\ell};\boldsymbol{u}_{1}^{\ell}\right)
\end{array}\right]\left[\begin{array}{c}
\boldsymbol{u}_{1}^{\ell+1}-\boldsymbol{u}_{1}^{\ell}\\
\\
\\
\boldsymbol{u}_{2}^{\ell+1}-\boldsymbol{u}_{2}^{\ell}
\end{array}\right]=-\left[\begin{array}{c}
\boldsymbol{f}_{1}\left(\boldsymbol{u}_{1}^{\ell};\boldsymbol{u}_{2}^{\ell}\right)\\
\\
\\
\boldsymbol{f}_{2}\left(\boldsymbol{u}_{2}^{\ell};\boldsymbol{u}_{1}^{\ell}\right)
\end{array}\right]
\end{equation}
to compute the solution updates. In general, developing a solver for
this linear system from scratch is quite challenging. Moreover, if
legacy solvers for each subproblem: $\boldsymbol{f}_{1}=\boldsymbol{0}$
and $\boldsymbol{f}_{2}=\boldsymbol{0}$, are available, extensive
modifications to their respective source codes would be required to 
construct the Jacobian in \hyperlink{eq22}{Eq. 2.2}.
Furthermore, the quadratic convergence rate of Newton's method is
not always guaranteed in practice. These limitations can be primarily
attributed to the bidirectional nature of the coupled PDE system and
the influence of the off-diagonal blocks ${\displaystyle \frac{\partial\boldsymbol{f}_{1}}{\partial\boldsymbol{u}_{2}}}$
and ${\displaystyle \frac{\partial\boldsymbol{f}_{2}}{\partial\boldsymbol{u}_{1}}}$
in the linear system. Variants of Newton's method such as Broyden's
method [\hyperlink{ref17}{17}], Gauss-Newton [\hyperlink{ref18}{18}] and Levenberg-Marquardt [\hyperlink{ref19}{19}],
would also be affected by these limitations, if implemented in this monolithic
fashion.

Therefore, modular approaches, for example, block-Gauss-Seidel (BGS)
[\hyperlink{ref20}{20}], can overcome some of these limitations of monolithic methods,
and are often preferred in practice. The primary advantage of modularization
is that legacy solvers with disparate discretization and solution
techniques can be easily coupled with minimal modifications to their
source codes. In general, the BGS method is an iterative method which
can be represented as 
\begin{equation}
\boldsymbol{u}_{1}^{\ell+1}=\boldsymbol{m}_{1}\left(\boldsymbol{u}_{1}^{\ell},\boldsymbol{u}_{2}^{\ell}\right),\ \boldsymbol{u}_{2}^{\ell+1}=\boldsymbol{m}_{2}\left(\boldsymbol{u}_{2}^{\ell},\boldsymbol{u}_{1}^{\ell+1}\right),
\end{equation}
where $\boldsymbol{m}_{1}\in\mathbb{R}^{n_{1}}$ and $\boldsymbol{m}_{2}\in\mathbb{R}^{n_{2}}$
denote the computational modules. Individually, these modules can
also be used to solve $\boldsymbol{f}_{1}=\boldsymbol{0}$, given
$\boldsymbol{u}_{2}$ and $\boldsymbol{f}_{2}=\boldsymbol{0}$, given
$\boldsymbol{u}_{1}$ respectively.

\subsection{Stochastic model}

We now consider the case where $\boldsymbol{f}_{1}$ and $\boldsymbol{f}_{2}$
are each dependent on a set of random input parameters, denoted as
$\boldsymbol{\xi}_{1}\in\Xi_{1}\subseteq\mathbb{R}^{s_{1}}$ and $\boldsymbol{\xi}_{2}\in\Xi_{2}\subseteq\mathbb{R}^{s_{2}}$
respectively. Let $\mu_{1}:\Xi_{1}\rightarrow\mathbb{R}^{+}$, $\mu_{2}:\Xi_{2}\rightarrow\mathbb{R}^{+}$
denote the respective probability density functions of $\boldsymbol{\xi}_{1}$ and $\boldsymbol{\xi}_{2}$. Moreover,
we assume that the input parameters as well as their individual elements
are statistically independent. Let $\Xi\equiv\Xi_{1}\times\Xi{}_{2}$
denote the combined parameter space with dimension $s=s_{1}+s_{2}$
and $\mu:\Xi\rightarrow\mathbb{R}^{+}:\forall\boldsymbol{\xi}\in\Xi$,
\begin{equation}
\mu\left(\boldsymbol{\xi}\right)=\mu_{1}\left(\boldsymbol{\xi}_{1}\left(\boldsymbol{\xi}\right)\right)\mu_{2}\left(\boldsymbol{\xi}_{2}\left(\boldsymbol{\xi}\right)\right)
\end{equation}
denote the corresponding joint probability density function.

Retaining the structure and dimensions of the deterministic algebraic
system in \hyperlink{eq21}{Eq. 2.1}, the stochastic nonlinear system of equations are
now formulated as follows. $\forall\boldsymbol{\xi}\in\Xi$ 
\hypertarget{eq25}{}
\begin{align}
 & \boldsymbol{f}_{1}\left(\boldsymbol{u}_{1}\left(\boldsymbol{\xi}\right);\boldsymbol{u}_{2}\left(\boldsymbol{\xi}\right),\boldsymbol{\xi}_{1}\left(\boldsymbol{\xi}\right)\right)=\boldsymbol{0},\nonumber \\
 & \boldsymbol{f}_{2}\left(\boldsymbol{u}_{2}\left(\boldsymbol{\xi}\right);\boldsymbol{u}_{1}\left(\boldsymbol{\xi}\right),\boldsymbol{\xi}_{2}\left(\boldsymbol{\xi}\right)\right)=\boldsymbol{0}.
\end{align}

We assume that $\forall\boldsymbol{\xi}\in\Xi$, the iterations
\begin{align}
 & \boldsymbol{u}_{1}^{\ell+1}\left(\boldsymbol{\xi}\right)=\boldsymbol{m}_{1}\left(\boldsymbol{u}_{1}^{\ell}\left(\boldsymbol{\xi}\right),\boldsymbol{u}_{2}^{\ell}\left(\boldsymbol{\xi}\right),\boldsymbol{\xi}_{1}\left(\boldsymbol{\xi}\right)\right),\nonumber \\
 & \boldsymbol{u}_{2}^{\ell+1}\left(\boldsymbol{\xi}\right)=\boldsymbol{m}_{2}\left(\boldsymbol{u}_{2}^{\ell}\left(\boldsymbol{\xi}\right),\boldsymbol{u}_{1}^{\ell+1}\left(\boldsymbol{\xi}\right),\boldsymbol{\xi}_{2}\left(\boldsymbol{\xi}\right)\right)
\end{align}
are to converge to a solution of \hyperlink{eq25}{Eq. 2.5}. 
A general formulation of a two-component coupled stochastic algebraic system \hypertarget{eq27}{}
\begin{align}
 & \boldsymbol{f}_{1}\left(\boldsymbol{u}_{1}\left(\boldsymbol{\xi}\right);\boldsymbol{v}_{2}\left(\boldsymbol{\xi}\right),\boldsymbol{\xi}_{1}\left(\boldsymbol{\xi}\right)\right)=\boldsymbol{0},\ \boldsymbol{v}_{1}\left(\boldsymbol{\xi}\right)=\boldsymbol{g}_{1}\left(\boldsymbol{u}_{1}\left(\boldsymbol{\xi}\right)\right),\nonumber \\
 & \boldsymbol{f}_{2}\left(\boldsymbol{u}_{2}\left(\boldsymbol{\xi}\right);\boldsymbol{v}_{1}\left(\boldsymbol{\xi}\right),\boldsymbol{\xi}_{2}\left(\boldsymbol{\xi}\right)\right)=\boldsymbol{0},\ \boldsymbol{v}_{2}\left(\boldsymbol{\xi}\right)=\boldsymbol{g}_{2}\left(\boldsymbol{u}_{2}\left(\boldsymbol{\xi}\right)\right),
\end{align}
where $\boldsymbol{g}_{1}\in\mathbb{R}^{m_{1}}$ and $\boldsymbol{g}_{2}\in\mathbb{R}^{m_{2}}$
denote the coupling or interface functions can be solved with a modified BGS method, wherein $\forall\boldsymbol{\xi}\in\Xi$,
the iterations 
\begin{align}
 \boldsymbol{u}_{1}^{\ell+1}\left(\boldsymbol{\xi}\right)=\boldsymbol{m}_{1}\left(\boldsymbol{u}_{1}^{\ell}\left(\boldsymbol{\xi}\right),\boldsymbol{v}_{2}^{\ell}\left(\boldsymbol{\xi}\right),\boldsymbol{\xi}_{1}\left(\boldsymbol{\xi}\right)\right),&\  \boldsymbol{v}_{1}^{\ell+1}=\boldsymbol{g}_{1}\left(\boldsymbol{u}_{1}^{\ell+1}\right),\nonumber \\
 \boldsymbol{u}_{2}^{\ell+1}\left(\boldsymbol{\xi}\right)=\boldsymbol{m}_{2}\left(\boldsymbol{u}_{2}^{\ell}\left(\boldsymbol{\xi}\right),\boldsymbol{v}_{1}^{\ell+1}\left(\boldsymbol{\xi}\right),\boldsymbol{\xi}_{2}\left(\boldsymbol{\xi}\right)\right),
&\  \boldsymbol{v}_{2}^{\ell+1}=\boldsymbol{g}_{2}\left(\boldsymbol{u}_{2}^{\ell+1}\right)
\end{align}
converge to a solution of \hyperlink{eq27}{Eq. 2.7}.

\subsection{Generalized polynomial chaos}

$\forall i\in\left\{ 1,2\right\} $, let $\boldsymbol{G}_{i}\in\mathbb{R}^{n_{i}\times n_{i}}$
denote the symmetric-positive-definite Gramian matrix [\hyperlink{ref21}{21}] corresponding
to $\boldsymbol{u}_{i}$. Moreover, let 
\begin{equation}
\mathcal{L}_{i}^{2}\left(\Xi\right)=\left\{ \boldsymbol{u}:\Xi\rightarrow\mathbb{R}^{n_{i}}:\int_{\Xi}\boldsymbol{u}\left(\boldsymbol{\xi}\right)^{\mathbf{T}}\boldsymbol{G}_{i}\boldsymbol{u}\left(\boldsymbol{\xi}\right)\mu\left(\boldsymbol{\xi}\right)d\boldsymbol{\xi}<\infty\right\} 
\end{equation}
denote the space of $\mu$-weighted, $\boldsymbol{G}_{i}$-square
integrable functions that map from $\Xi$ to $\mathbb{R}^{n_{i}}$.
If $\boldsymbol{u}_{i}\in\mathcal{L}_{i}^{2}\left(\Xi\right)$, then it can be represented exactly
as an infinite polynomial series as follows. $\forall\boldsymbol{\xi}\in\Xi$,
\hypertarget{eq210}{}
\begin{equation}
\boldsymbol{u}_{i}\left(\boldsymbol{\xi}\right)=\sum_{j\geq0}\hat{\boldsymbol{u}}_{i,j}\psi_{j}\left(\boldsymbol{\xi}\right),
\end{equation}
where $\left\{ \psi_{j}:\Xi\rightarrow\mathbb{R}\right\} _{j\geq0}$
denotes the set of $\mu$-orthonormal polynomials, and a basis for
all $ \mu$-weighted, square-integrable scalar functions in $\Xi$.
The statistical independence of the coordinate directions in $\Xi$
implies that each basis polynomial is a product of $s$ univariate
orthonormal polynomials [\hyperlink{ref22}{22}], which in turn can be precomputed
using the Golub-Welsch algorithm [\hyperlink{ref23}{23}]. Moreover, the indexing
in the polynomial series is assumed to follow a total degree ordering, such that
\begin{equation}
\mathrm{deg}\left(\psi_{j}\right)\geq\mathrm{deg}\left(\psi_{k}\right)\Leftrightarrow j\geq k\geq0,
\end{equation}
where $\forall j\geq0$, $\mathrm{deg}\left(\psi_{j}\right)$ denotes
the total degree of $\psi_{j}$. 

Consequently, a gPC [\hyperlink{ref8}{8}] approximation $\boldsymbol{u}_{i}^{p}\approx\boldsymbol{u}_{i}:\forall\boldsymbol{\xi}\in\Xi$,
\begin{equation}
\boldsymbol{u}_{i}^{p}\left(\boldsymbol{\xi}\right)=\sum_{j=0}^{P}\hat{\boldsymbol{u}}_{i,j}\psi_{j}\left(\boldsymbol{\xi}\right)=\hat{\boldsymbol{U}}_{i}\boldsymbol{\psi}\left(\boldsymbol{\xi}\right)
\end{equation}
with order $p\geq0$ and $P+1$ gPC coefficients can be formulated
as a finite truncation of \hyperlink{eq210}{Eq. 2.10}, where $\hat{\boldsymbol{U}}_{i}\equiv\hat{\boldsymbol{U}}_{i}^{p}=\left[\begin{array}{ccc}
\hat{\boldsymbol{u}}_{i,0} & \cdots & \hat{\boldsymbol{u}}_{i,P}\end{array}\right]\in\mathbb{R}^{n_{i}\times\left(P+1\right)}$ denotes the gPC coefficient matrix and $\boldsymbol{\psi}\equiv\boldsymbol{\psi}^{p}=\left[\begin{array}{ccc}
\psi_{0} & \cdots & \psi_{P}\end{array}\right]^{\mathbf{T}}:\Xi\rightarrow\mathbb{R}^{P+1}$ denotes the basis vector. If the truncation is isotropic
and based on the total degree, then $p$ and $P$ are related as follows.
\begin{equation}
P+1=\frac{\left(p+s\right)!}{p!s!}.
\end{equation}

The orthonormality of the
elements in $\boldsymbol{\psi}$ yields the spectral
projection formula 
\hypertarget{eq214}{}
\begin{equation}
\hat{\boldsymbol{U}}_{i}=\int_{\Xi}\boldsymbol{u}_{i}\left(\boldsymbol{\xi}\right)\boldsymbol{\psi}\left(\boldsymbol{\xi}\right)^{\mathbf{T}}\mu\left(\boldsymbol{\xi}\right)d\boldsymbol{\xi}.
\end{equation}

Moreover, the Cameron-Martin theorem [\hyperlink{ref24}{24}] states that if $\boldsymbol{u}_{i}$
is infinitely regular in $\Xi$, then as $p\rightarrow\infty$, the
gPC approximation $\boldsymbol{u}_{i}^{p}$ converges
exponentially to $\boldsymbol{u}_{i}$, in the mean-square sense.
To state this formally, $\exists\chi^{*}>0,\rho^{*}>1$, such that $\forall p\geq 0$,
\begin{equation}
\sqrt{\int_{\Xi}\left\Vert \boldsymbol{u}_{i}^{p}\left(\boldsymbol{\xi}\right)-\boldsymbol{u}_{i}\left(\boldsymbol{\xi}\right)\right\Vert _{\boldsymbol{G}_{i}}^{2}\mu\left(\boldsymbol{\xi}\right)d\boldsymbol{\xi}}\leq\chi^{*}\rho^{-p}
\end{equation}
for some $\rho>\rho^{*}$. Furthermore, a corollary to the theorem states
that if the regularity of $\boldsymbol{u}_{i}$ is $k$, then the
asymptotic rate of convergence is polynomial and the upper bound in
the mean-square error would be $\mathcal{O}\left(p^{-k}\right)$. 

When compared to repeated executions of a costly numerical PDE solver, the gPC approximation $\boldsymbol{u}_{i}^{p}$
can be used as a significantly cheaper surrogate model for computing
statistical quantities of interest such as the moments and probability
density functions of $\boldsymbol{u}_{i}$ via exhaustive sampling.
In particular, approximations of the first two moments can be computed
directly from the gPC coefficients as follows.
\begin{align}
 & \mathbb{E}\left(\boldsymbol{u}_{i}\right)\approx\mathbb{E}\left(\boldsymbol{u}_{i}^{p}\right)=\hat{\boldsymbol{u}}_{i,0}, \nonumber\\
 & \mathrm{Cov}\left(\boldsymbol{u}_{i},\boldsymbol{u}_{i}\right)\approx\mathrm{Cov}\left(\boldsymbol{u}_{i}^{p},\boldsymbol{u}_{i}^{p}\right)=\hat{\boldsymbol{U}}_{i}\hat{\boldsymbol{U}}_{i}^{\mathbf{T}}-\hat{\boldsymbol{u}}_{i,0}\hat{\boldsymbol{u}}_{i,0}^{\mathbf{T}}.
\end{align}

Moreover, the probability density function of any related quantity of interest can be computed
with the kernel density estimation (KDE) method [\hyperlink{ref25}{25}], by generating a large number of samples of the cheaper polynomial surrogates. Therefore, prior to
any uncertainty analysis, the gPC coefficient matrices of the solutions
must be computed via uncertainty propagation. In this work, we focus on the non-intrusive projection method [\hyperlink{ref26}{26}], which enables
the reuse of the solver components $\boldsymbol{m}_{1}$ and $\boldsymbol{m}_{2}$.
\subsection{Non-intrusive spectral projection}
Non-intrusive spectral projection (NISP), is a  method to approximate the gPC coefficients
by approximating the integration in the spectral projection formula
in \hyperlink{eq214}{Eq. 2.14} using a quadrature rule in in $\Xi$, denoted as $\left\{ \left(\boldsymbol{\xi}^{\left( j\right)},w^{\left( j\right)}\right)\right\} _{j=1}^{Q}$.
Therefore, $\forall i\in\left\{ 1,2\right\} $, $\hat{\boldsymbol{U}}_{i}:$
\begin{equation}
\hat{\boldsymbol{U}}_{i}\approx\sum_{j=1}^{Q}w^{\left( j\right)}\boldsymbol{u}_{i}\left(\boldsymbol{\xi}^{\left( j\right)}\right)\boldsymbol{\psi}\left(\boldsymbol{\xi}\right)^{\mathbf{T}},
\end{equation}
If $q\geq0$ denotes the level of the quadrature rule, then all polynomials with
a total degree $\leq2q+1$ can be numerically integrated using the quadrature rule, up to machine precision. The relation between $Q$ and $q$ depends
on the type of quadrature rule implemented. For instance, employing a full grid
tensorization of univariate Gauss-quadrature rules yields the relation
$Q=\left(q+1\right)^{s}$, while for the same level $q$, sparse grids [\hyperlink{ref27}{27}] exhibit a much
slower, albeit exponential growth in $Q$. This exponential growth phenomenon
is commonly known as the curse of dimension.

\subsubsection{Algorithm and computational cost}

\hyperlink{alg1}{Algorithm 1} describes the standard NISP based uncertainty
propagation method.

$\forall i\in\left\{ 1,2\right\} $, the precomputed matrix $\hat{\boldsymbol{\Xi}}_{i}\equiv\hat{\boldsymbol{\Xi}}_{i}^{p}\in\mathbb{R}^{s_{i}\times\left(P+1\right)}:$

\begin{align}
\hat{\boldsymbol{\Xi}}_{i} & =\int_{\Xi}\boldsymbol{\xi}_{i}\left(\boldsymbol{\xi}\right)\boldsymbol{\psi}\left(\boldsymbol{\xi}\right)^{\mathbf{T}}\mu\left(\boldsymbol{\xi}\right)d\boldsymbol{\xi}
=\sum_{j=1}^{Q}w^{\left( j\right)}\boldsymbol{\xi}_{i}\left(\boldsymbol{\xi}^{\left( j\right)}\right)\boldsymbol{\psi}\left(\boldsymbol{\xi}^{\left( j\right)}\right)^{\mathbf{T}}
\end{align}
denotes the gPC coefficient matrix of the random input $\boldsymbol{\xi}_{i}$. 

We assume that the computational costs are dominated by the repeated execution of module 
operators $\boldsymbol{m}_{1}$ and $\boldsymbol{m}_{2}$. Let $\bar{\mathcal{C}}_{1}$
denote the average cost of executing of $\boldsymbol{m}_{1}$
and $\bar{\mathcal{C}}_{2}$ denote the average cost of executing of
$\boldsymbol{m}_{2}$. Therefore, the computational cost of the standard
NISP method 
\begin{equation}
C_{s}\approx\mathcal{O}\left(\left(\bar{\mathcal{C}}_{1}+\bar{\mathcal{C}}_{2}\right)Q\right)
\end{equation}
would grow exponentially with respect to the dimension
$s$ and order $p$. To mitigate these costs, we propose a reduced NISP based uncertainty
propagation method, which will be described in the next section.
\RestyleAlgo{boxruled}
\begin{algorithm}
\hypertarget{alg1}{}
\caption{Standard NISP based uncertainty propagation for a two-module multi-physics system}
\SetAlgoLined

\SetKwInOut{Input}{inputs}\SetKwInOut{Output}{outputs}
\DontPrintSemicolon
\Input{$\mu_{1}$, $\mu_{2}$, order $p\geq0$, level $q\geq p$, $\hat{\boldsymbol{U}}_{1}^{0}$, $\hat{\boldsymbol{U}}_{2}^{0}$}
\Output{$\hat{\boldsymbol{U}}_{1}$, $\hat{\boldsymbol{U}}_{2}$}
\textbf{precompute}:: $\left\{ \left(w^{\left( j\right)},\boldsymbol{\psi}\left(\boldsymbol{\xi}^{\left( j\right)}\right)\right)\right\} _{j=1}^{Q}$,
$\hat{\boldsymbol{\Xi}}_{1}$, $\hat{\boldsymbol{\Xi}}_{2}$

\textbf{$\ell\leftarrow0$}

\Repeat{$\hat{\boldsymbol{U}}_{1}^{\ell}$, $\hat{\boldsymbol{U}}_{2}^{\ell}$ $\mathrm{converge}$}{

\textbf{}$\hat{\boldsymbol{U}}_{1}^{\ell+1}\leftarrow\boldsymbol{0}$

\For{$ j\leftarrow 1$ \KwTo $Q$}{

\textbf{$\boldsymbol{u}_{1}\leftarrow\boldsymbol{m}_{1}\left(\hat{\boldsymbol{U}}_{1}^{\ell}\boldsymbol{\psi}\left(\boldsymbol{\xi}^{\left( j\right)}\right),\hat{\boldsymbol{U}}_{2}^{\ell}\boldsymbol{\psi}\left(\boldsymbol{\xi}^{\left( j\right)}\right),\hat{\boldsymbol{\Xi}}_{1}\boldsymbol{\psi}\left(\boldsymbol{\xi}^{\left( j\right)}\right)\right)$}

\textbf{$\hat{\boldsymbol{U}}_{1}^{\ell+1}\leftarrow\hat{\boldsymbol{U}}_{1}^{\ell+1}+w^{\left( j\right)}\boldsymbol{u}_{1}\boldsymbol{\psi}\left(\boldsymbol{\xi}^{\left( j\right)}\right)^{\mathbf{T}}$}

}

\textbf{}$\hat{\boldsymbol{U}}_{2}^{\ell+1}\leftarrow\boldsymbol{0}$ 

\For{$ j\leftarrow 1$ \KwTo $Q$}{

\textbf{$\boldsymbol{u}_{2}\leftarrow\boldsymbol{m}_{2}\left(\hat{\boldsymbol{U}}_{2}^{\ell}\boldsymbol{\psi}\left(\boldsymbol{\xi}^{\left( j\right)}\right),\hat{\boldsymbol{U}}_{1}^{\ell+1}\boldsymbol{\psi}\left(\boldsymbol{\xi}^{\left( j\right)}\right),\hat{\boldsymbol{\Xi}}_{2}\boldsymbol{\psi}\left(\boldsymbol{\xi}^{\left( j\right)}\right)\right)$};

\textbf{$\hat{\boldsymbol{U}}_{2}^{\ell+1}\leftarrow\hat{\boldsymbol{U}}_{2}^{\ell+1}+w^{\left( j\right)}\boldsymbol{u}_{2}\boldsymbol{\psi}\left(\boldsymbol{\xi}^{\left( j\right)}\right)^{\mathbf{T}}$}

}

\textbf{$\ell\leftarrow \ell+1$}

}

\end{algorithm}

\section{Reduced NISP based uncertainty propagation}
\hypertarget{sec3}{}
Our proposed method is a modification of \hyperlink{alg1}{Algorithm 1} with the addition
of two intermediate computational steps at each iteration: (1) a dimension reduction
routine,  and (2) an order reduction routine.

\subsection{Dimension reduction}
\hypertarget{sec31}{}
At each iteration $\ell$, let $\boldsymbol{y}_{1}^{\ell}\equiv\left[\boldsymbol{u}_{1}^{\ell};\boldsymbol{u}_{2}^{\ell};\boldsymbol{\xi}_{1}\right]:\Xi\rightarrow\mathbb{R}^{r_{1}}$
and $\boldsymbol{y}_{2}^{\ell}\equiv\left[\boldsymbol{u}_{2}^{\ell};\boldsymbol{u}_{1}^{\ell+1};\boldsymbol{\xi}_{2}\right]:\Xi\rightarrow\mathbb{R}^{r_{2}}$
denote the input data that enter $\boldsymbol{m}_{1}$ and $\boldsymbol{m}_{2}$
respectively. Moreover, let $\boldsymbol{\Gamma}_{1},\boldsymbol{\Gamma}_{2}:$
\begin{equation}
\boldsymbol{\Gamma}_{1}=\left[\begin{array}{ccc}
\boldsymbol{G}_{1} & \boldsymbol{0} & \boldsymbol{0}\\
\boldsymbol{0} & \boldsymbol{G}_{2} & \boldsymbol{0}\\
\boldsymbol{0} & \boldsymbol{0} & \boldsymbol{I}_{s_{1}}
\end{array}\right]\in\mathbb{R}^{r_{1}\times r_{1}},\boldsymbol{\Gamma}_{2}=\left[\begin{array}{ccc}
\boldsymbol{G}_{2} & \boldsymbol{0} & \boldsymbol{0}\\
\boldsymbol{0} & \boldsymbol{G}_{1} & \boldsymbol{0}\\
\boldsymbol{0} & \boldsymbol{0} & \boldsymbol{I}_{s_{2}}
\end{array}\right]\in\mathbb{R}^{r_{2}\times r_{2}}
\end{equation}
denote the Gramian matrices and 
\begin{equation}
\hat{\boldsymbol{Y}}_{1}^{\ell}\equiv\hat{\boldsymbol{Y}}_{1}^{\ell,p}=\left[\begin{array}{c}
\hat{\boldsymbol{U}}_{1}^{\ell}\\
\hat{\boldsymbol{U}}_{2}^{\ell}\\
\hat{\boldsymbol{\Xi}}_{1}
\end{array}\right]\in\mathbb{R}^{r_{1}\times\left(P+1\right)},\hat{\boldsymbol{Y}}_{2}^{\ell}\equiv\hat{\boldsymbol{Y}}_{2}^{\ell,p}=\left[\begin{array}{c}
\hat{\boldsymbol{U}}_{2}^{\ell}\\
\hat{\boldsymbol{U}}_{1}^{\ell+1}\\
\hat{\boldsymbol{\Xi}}_{2}
\end{array}\right]\in\mathbb{R}^{r_{2}\times\left(P+1\right)}.
\end{equation}
as the gPC coefficient matrices corresponding to $\boldsymbol{y}_{1}^{\ell}$
and $\boldsymbol{y}_{2}^{\ell}$ respectively. Subsequently,
$\forall i\in\left\{ 1,2\right\} $, we construct $\tilde{\boldsymbol{Y}}_{i}^{\ell}:$ 
\begin{equation}
\tilde{\boldsymbol{Y}}_{i}^{\ell}=\boldsymbol{\Gamma}_{i}^{\frac{1}{2}}\left[\begin{array}{ccc}
\hat{\boldsymbol{y}}_{i,1}^{\ell} & \cdots & \hat{\boldsymbol{y}}_{i,P}^{\ell}\end{array}\right],
\end{equation}
by deleting the first column of
$\hat{\boldsymbol{Y}}_{i}^{\ell}$, and left multiplying
the resultant matrix with $\boldsymbol{\Gamma}_{i}^{\frac{1}{2}}$. The singular value decomposition (SVD) of $\tilde{\boldsymbol{Y}}_{i}^{\ell}:$ 

\begin{equation}
\tilde{\boldsymbol{Y}}_{i}^{\ell}=\tilde{\boldsymbol{\Upsilon}}_{i}^{\ell}\boldsymbol{\Sigma}_{i}^{\ell}\left(\tilde{\boldsymbol{\Theta}}_{i}^{\ell}\right)^{\mathbf{T}}
\end{equation}
is computed, wherein the constituent matrices have the following
structure and dimensions.

\begin{align}
\tilde{\boldsymbol{\Upsilon}}_{i}^{\ell} & =\left[\begin{array}{ccc}
\tilde{\boldsymbol{\upsilon}}_{i,1}^{\ell} & \cdots & \tilde{\boldsymbol{\upsilon}}_{i,\min\left\{ P,r_{i}\right\} }^{\ell}\end{array}\right]\in\mathbb{R}^{r_{i}\times\min\left\{ P,r_{i}\right\} },\\
\boldsymbol{\Sigma}_{i}^{\ell} & =\left[\begin{array}{ccc}
\sigma_{i,1}^{\ell} & 0 & 0\\
0 & \ddots & 0\\
0 & 0 & \sigma_{i,\min\left\{ P,r_{i}\right\} }^{\ell}
\end{array}\right]\in\mathbb{R}^{\min\left\{ P,r_{i}\right\} \times\min\left\{ P,r_{i}\right\} },\\
\tilde{\boldsymbol{\Theta}}_{i}^{\ell} & =\left[\begin{array}{ccc}
\hat{\theta}_{i,1,1}^{\ell} & \cdots & \hat{\theta}_{i,1,\min\left\{ P,r_{i}\right\} }^{\ell}\\
\vdots &  & \vdots\\
\hat{\theta}_{i,P,1}^{\ell} & \cdots & \hat{\theta}_{i,P,\min\left\{ P,r_{i}\right\} }^{\ell}
\end{array}\right]\in\mathbb{R}^{P\times\min\left\{ P,r_{i}\right\} }.
\end{align}

We then construct 
\begin{equation}
\boldsymbol{\Upsilon}_{i}^{\ell}=\boldsymbol{\Gamma}_{i}^{-\frac{1}{2}}\tilde{\boldsymbol{\Upsilon}}_{i}^{\ell}=\left[\begin{array}{ccc}
\boldsymbol{\upsilon}_{i,1}^{\ell} & \cdots & \boldsymbol{\upsilon}_{i,\min\left\{ P,r_{i}\right\} }^{\ell}\end{array}\right],
\end{equation}
and 
\begin{equation}
\hat{\boldsymbol{\Theta}}_{i}^{\ell}=\left[\begin{array}{ccc}
0 & \cdots & 0\\
 & \tilde{\boldsymbol{\Theta}}_{i}^{\ell}
\end{array}\right]=\left[\begin{array}{ccc}
\hat{\boldsymbol{\theta}}_{i,1}^{\ell} \\
\vdots \\
\hat{\boldsymbol{\theta}}_{i,\min\left\{ P,r_{i}\right\} }^{\ell} 
\end{array}\right]^{\mathbf{T}}.
\end{equation}

Subsequently, the gPC approximation $\boldsymbol{y}_{i}^{\ell,p}$
can be rewritten as follows. $\forall\boldsymbol{\xi}\in\Xi$,
\hypertarget{eq310}{}
\begin{equation}
\boldsymbol{y}_{i}^{\ell,p}\left(\boldsymbol{\xi}\right)=\hat{\boldsymbol{y}}_{i,0}^{\ell}+\sum_{j=1}^{\min\left\{ P,r_{i}\right\} }\sigma_{i,j}^{\ell}\boldsymbol{\upsilon}_{i,j}^{\ell}\theta_{i,j}^{\ell}\left(\boldsymbol{\xi}\right),
\end{equation}
where $\forall1\leq j\leq\min\left\{ P,r_{i}\right\} $, 
\begin{align}
\theta_{i,j}^{\ell}\left(\boldsymbol{\xi}\right) & =\hat{\boldsymbol{\theta}}_{i,j}^{\ell}\boldsymbol{\psi}\left(\boldsymbol{\xi}\right)
\end{align}

denotes the $j$-th uncorrelated random variable with zero mean and unit variance that enters the component  $\boldsymbol{m}_{i}$. The expansion in \hyperlink{eq310}{Eq. 3.10} is the finite-dimensional variant of the Karhunen-Loeve (KL) expansion [\hyperlink{ref28}{28}], which is widely used
as an approximation of spatiotemporal random fields.

\subsubsection{Reduced dimensional KL expansion}

$\forall i\in\left\{ 1,2\right\} $ and iteration $\ell$, the expansion
(3.10) can be truncated by retaining $d_{i}\equiv d_{i}^{\ell}$
terms to define an approximation $\boldsymbol{y}_{i}^{\ell,d_{i}}\equiv\boldsymbol{y}_{i}^{\ell,p,d_{i}}:\forall\boldsymbol{\xi}\in\Xi$, 
\begin{align}
\boldsymbol{y}_{i}^{\ell,p}\left(\boldsymbol{\xi}\right)\approx\boldsymbol{y}_{i}^{\ell,d_{i}}\left(\boldsymbol{\xi}\right) & =\hat{\boldsymbol{y}}_{i,0}^{\ell}+\sum_{j=1}^{d_{i}}\sigma_{i,j}^{\ell}\boldsymbol{\upsilon}_{i,j}^{\ell}\theta_{i,j}^{\ell}\left(\boldsymbol{\xi}\right).
\end{align}

Let $\boldsymbol{\theta}_{i}^{\ell}\equiv\left[\begin{array}{ccc}
\theta_{i,1}^{\ell} & \cdots & \theta_{i,d_{i}}^{\ell}\end{array}\right]^{\mathbf{T}}:\Xi\rightarrow\Theta_{i}^{\ell}\subseteq\mathbb{R}^{d_{i}}$ denote the reduced dimensional random vector with a  probability density function $\nu_{i}^{\ell}:\Theta_{i}^{\ell}\rightarrow\mathbb{R}^{+}$.
Moreover, we let $\boldsymbol{z}_{i}^{\ell}:\Theta_{i}^{\ell}\rightarrow\mathbb{R}^{r_{i}}$
denote the affine map which defines the composition $\boldsymbol{y}_{i}^{\ell,d_{i}}=\boldsymbol{z}_{i}^{\ell}\circ\boldsymbol{\theta}_{i}^{\ell}$, such that $\forall\boldsymbol{\theta}\in\Theta_{i}^{\ell}$,
\begin{equation}
\boldsymbol{z}_{i}^{\ell}\left(\boldsymbol{\theta}\right)=\bar{\boldsymbol{z}}_{i}^{\ell}+\tilde{\boldsymbol{Z}}_{i}^{\ell}\boldsymbol{\theta},
\end{equation}
where 
\begin{equation}
\bar{\boldsymbol{z}}_{i}^{\ell}=\hat{\boldsymbol{y}}_{i,0}^{\ell}\in\mathbb{R}^{r_{i}}
\end{equation}
and
\begin{equation}
\tilde{\boldsymbol{Z}}_{i}^{\ell}=\left[\begin{array}{ccc}
\boldsymbol{\upsilon}_{i,1}^{\ell} & \cdots & \boldsymbol{\upsilon}_{i,d_{i}}^{\ell}\end{array}\right]\left[\begin{array}{ccc}
\sigma_{i,1}^{\ell}\\
 & \ddots\\
 &  & \sigma_{i,d_{i}}^{\ell}
\end{array}\right]\in\mathbb{R}^{r_{i}\times d_{i}}.
\end{equation}

Furthermore, we extract the following subvector and submatrix blocks
from $\bar{\boldsymbol{z}}_{i}^{\ell}$ and $\tilde{\boldsymbol{Z}}_{i}^{\ell}$
respectively. 
\begin{equation}
\bar{\boldsymbol{z}}_{1}^{\ell}=\left[\begin{array}{c}
\bar{\boldsymbol{u}}_{1,1}^{\ell}\\
\bar{\boldsymbol{u}}_{2,1}^{\ell}\\
\bar{\boldsymbol{\xi}}_{1,1}^{\ell}
\end{array}\right],\bar{\boldsymbol{z}}_{2}^{\ell}=\left[\begin{array}{c}
\bar{\boldsymbol{u}}_{2,2}^{\ell}\\
\bar{\boldsymbol{u}}_{1,2}^{\ell}\\
\bar{\boldsymbol{\xi}}_{2,2}^{\ell}
\end{array}\right],
\end{equation}
 
\begin{equation}
\tilde{\boldsymbol{Z}}_{1}^{\ell}=\left[\begin{array}{c}
\tilde{\boldsymbol{U}}_{1,1}^{\ell}\\
\tilde{\boldsymbol{U}}_{2,1}^{\ell}\\
\tilde{\boldsymbol{\Xi}}_{1,1}^{\ell}
\end{array}\right],\tilde{\boldsymbol{Z}}_{2}^{\ell}=\left[\begin{array}{c}
\tilde{\boldsymbol{U}}_{2,2}^{\ell}\\
\tilde{\boldsymbol{U}}_{1,2}^{\ell}\\
\tilde{\boldsymbol{\Xi}}_{2,2}^{\ell}
\end{array}\right],
\end{equation}
where $\bar{\boldsymbol{u}}_{1,1}^{\ell},\bar{\boldsymbol{u}}_{1,2}^{\ell}\in\mathbb{R}^{n_{1}}$,
$\bar{\boldsymbol{u}}_{2,1}^{\ell},\bar{\boldsymbol{u}}_{2,2}^{\ell}\in\mathbb{R}^{n_{2}}$,
$\bar{\boldsymbol{\xi}}_{1,1}^{\ell}\in\mathbb{R}^{s_{1}}$,
$\bar{\boldsymbol{\xi}}_{2,2}^{\ell}\in\mathbb{R}^{s_{2}}$ and 
$\tilde{\boldsymbol{U}}_{1,1}^{\ell}\in\mathbb{R}^{n_{1}\times d_{1}}$,
$\tilde{\boldsymbol{U}}_{1,2}^{\ell}\in\mathbb{R}^{n_{1}\times d_{2}}$,
$\tilde{\boldsymbol{U}}_{2,1}^{\ell}\in\mathbb{R}^{n_{2}\times d_{1}}$,
$\tilde{\boldsymbol{U}}_{2,2}^{\ell}\in\mathbb{R}^{n_{2}\times d_{2}}$,
$\tilde{\boldsymbol{\Xi}}_{1,1}^{\ell}\in\mathbb{R}^{s_{1}\times d_{1}}$,
$\tilde{\boldsymbol{\Xi}}_{2,2}^{\ell}\in\mathbb{R}^{s_{2}\times d_{2}}$ define the affine maps corresponding to 
$\boldsymbol{u}_{1}, \boldsymbol{u}_{2}, \boldsymbol{\xi}_{1}, \boldsymbol{\xi}_{2}$ respectively.

\subsubsection{Selecting the reduced dimensions}

$\forall i\in\left\{ 1,2\right\} $, we prescribe a tolerance $\epsilon_{i,\mathrm{dim}}>0$
such that, at each iteration $\ell$, $d_{i}$ is selected as the minimum $k\in\mathbb{N}$ which satisfies 
\hypertarget{eq318}{}
\begin{equation}
\sqrt{\sum_{j=k+1}^{\min\left\{ P,r_{i}\right\} }\left(\sigma_{i,j}^{\ell}\right)^{2}}\leq\epsilon_{i,\mathrm{dim}}\sqrt{\sum_{j=1}^{\min\left\{ P,r_{i}\right\} }\left(\sigma_{i,j}^{\ell}\right)^{2}}.
\end{equation}

If $d_{i}<s$, then a reduced dimensional approximation of the input
data in module $\boldsymbol{m}_{i}$ exists, with an approximation
error $\leq\mathcal{O}\left(\epsilon_{i,\mathrm{dim}}\right)$.
\hyperlink{thm1}{Theorem 1} proves this error bound.
\subsubsection*{Theorem 1:}
\hypertarget{thm1}{}
\emph{$\forall i\in\left\{ 1,2\right\} $ and iteration $\ell$, the approximation
$\boldsymbol{y}_{i}^{\ell,d_{i}}$ satisfies the
inequality
\hypertarget{eq319}{}
\begin{equation}
\sqrt{\int_{\Xi}\left\Vert \boldsymbol{y}_{i}^{\ell}\mathbf{\left(\boldsymbol{\xi}\right)}-\boldsymbol{y}_{i}^{\ell,d_{i}}\left(\boldsymbol{\xi}\right)\right\Vert _{\boldsymbol{\Gamma}_{i}}^{2}\mu\left(\boldsymbol{\xi}\right)d\boldsymbol{\xi}}\leq\epsilon_{i,\mathrm{dim}}\sqrt{\int_{\Xi}\left\Vert \boldsymbol{y}_{i}^{\ell}\left(\boldsymbol{\xi}\right)-\hat{\boldsymbol{y}}_{i,0}^{\ell}\right\Vert _{\boldsymbol{\Gamma}_{i}}^{2}\mu\left(\boldsymbol{\xi}\right)d\boldsymbol{\xi}}.
\end{equation}
}

\subsubsection*{Proof:}
The square of the left hand side expression in \hyperlink{eq319}{Eq. 3.19} can be written
as follows. 
\hypertarget{eq320}{}
\begin{align}
 & \int_{\Xi}\left\Vert \boldsymbol{y}_{i}^{\ell}\mathbf{\left(\boldsymbol{\xi}\right)}-\boldsymbol{y}_{i}^{\ell,d_{i}}\left(\boldsymbol{\xi}\right)\right\Vert _{\boldsymbol{\Gamma}_{i}}^{2}\mu\left(\boldsymbol{\xi}\right)d\boldsymbol{\xi}=\int_{\Xi}\left\Vert \left(\sum_{j=d_{i}+1}^{\min\left\{ P,r_{i}\right\} }\sigma_{i,j}^{\ell}\boldsymbol{\upsilon}_{i,j}^{\ell}\theta_{i,j}^{\ell}\left(\boldsymbol{\xi}\right)\right)\right\Vert _{\boldsymbol{\Gamma}_{i}}^{2}\mu\left(\boldsymbol{\xi}\right)d\boldsymbol{\xi} \nonumber\\
 & =\int_{\Xi}\left(\left(\boldsymbol{\psi}\left(\boldsymbol{\xi}\right)^{\mathbf{T}}\left[\begin{array}{ccc}
\hat{\boldsymbol{\theta}}_{i,d_{i}+1}^{\ell} \\
\vdots \\
\hat{\boldsymbol{\theta}}_{i,\min\left\{ P,r_{i}\right\} }^{\ell} 
\end{array}\right]^{\mathbf{T}}\left[\begin{array}{ccc}
\sigma_{i,d_{i}+1}^{\ell} \\
 & \ddots\\
 &  & \sigma_{i,\min\left\{ P,r_{i}\right\} }^{\ell}
\end{array}\right] \right.\right. \nonumber\\
 & \left.\ \ \times\left[\begin{array}{ccc}
\tilde{\boldsymbol{\upsilon}}_{i,d_{i}+1}^{\ell} & \cdots & \tilde{\boldsymbol{\upsilon}}_{i,\min\left\{ P,r_{i}\right\} }^{\ell}\end{array}\right]^{\mathbf{T}}\boldsymbol{\Gamma}_{i}^{-\frac{1}{2}}\right)\boldsymbol{\Gamma}_{i} \left(\boldsymbol{\Gamma}_{i}^{-\frac{1}{2}}\left[\begin{array}{ccc}
\tilde{\boldsymbol{\upsilon}}_{i,d_{i}+1}^{\ell} & \cdots & \tilde{\boldsymbol{\upsilon}}_{i,\min\left\{ P,r_{i}\right\} }^{\ell}\end{array}\right]\right. \nonumber\\
 & \ \  \times\left.\left.\left[\begin{array}{ccc}
\sigma_{i,d_{i}+1}^{\ell}\\
 & \ddots\\
 &  & \sigma_{i,\min\left\{ P,r_{i}\right\} }^{\ell}
\end{array}\right]\left[\begin{array}{ccc}
\hat{\boldsymbol{\theta}}_{i,1}^{\ell} \nonumber\\
\vdots \\
\hat{\boldsymbol{\theta}}_{i,\min\left\{ P,r_{i}\right\} }^{\ell} 
\end{array}\right]\boldsymbol{\psi}\left(\boldsymbol{\xi}\right)\right)\right)\mu\left(\boldsymbol{\xi}\right)d\boldsymbol{\xi}\nonumber \\
 & =\mathrm{trace}\left(\left(\int_{\Xi}\boldsymbol{\psi}\left(\boldsymbol{\xi}\right)\boldsymbol{\psi}\left(\boldsymbol{\xi}\right)^{\mathbf{T}}\mu\left(\boldsymbol{\xi}\right)d\boldsymbol{\xi}\right)\left(\left[\begin{array}{ccc}
\hat{\boldsymbol{\theta}}_{i,d_{i}+1}^{\ell} \\
\vdots \\
\hat{\boldsymbol{\theta}}_{i,\min\left\{ P,r_{i}\right\} }^{\ell} 
\end{array}\right]^{\mathbf{T}}\left[\begin{array}{ccc}
\sigma_{i,d_{i}+1}^{\ell}\\
 & \ddots\\
 &  & \sigma_{i,\min\left\{ P,r_{i}\right\} }^{\ell}
\end{array}\right]^{2}\right.\right. \nonumber \\
 & \left.\left.\ \ \times\left[\begin{array}{ccc}
\hat{\boldsymbol{\theta}}_{i,d_{i}+1}^{\ell} \\
\vdots \\
\hat{\boldsymbol{\theta}}_{i,\min\left\{ P,r_{i}\right\} }^{\ell} 
\end{array}\right]\right)\right) =\mathrm{trace}\left(\left[\begin{array}{ccc}
\left(\sigma_{i,d_{i}+1}^{\ell}\right)^{2}\\
 & \ddots\\
 &  & \left(\sigma_{i,\min\left\{ P,r_{i}\right\} }^{\ell}\right)^{2}
\end{array}\right]\right)=\sum_{j=d_{i}+1}^{\min\left\{ P,r_{i}\right\} }\left(\sigma_{i,j}^{\ell}\right)^{2}.
\end{align}

Similarly, for the right hand side expression, we can show that
\hypertarget{eq321}{}
\begin{equation}
\int_{\Xi}\left\Vert \boldsymbol{y}_{i}^{\ell}\left(\boldsymbol{\xi}\right)-\hat{\boldsymbol{y}}_{i,0}^{\ell}\right\Vert _{\boldsymbol{\Gamma}_{i}}^{2}\mu\left(\boldsymbol{\xi}\right)d\boldsymbol{\xi}=\sum_{j=1}^{\min\left\{ P,r_{i}\right\} }\left(\sigma_{i,j}^{\ell}\right)^{2}.
\end{equation}

By substituting \hyperlink{eq320}{Eq. 3.20} and \hyperlink{eq321}{Eq. 3.21} into \hyperlink{eq318}{Eq. 3.18}, we arrive at \hyperlink{eq319}{Eq. 3.19}. $\square$

In the context of the modified coupled model with interface functions
in \hyperlink{eq27}{Eq. 2.7}, the dimension reduction procedure would be exactly the
same as described here with the exception that the input data are
formulated as as $\boldsymbol{y}_{1}^{\ell}\equiv\left[\boldsymbol{u}_{1}^{\ell};\boldsymbol{v}_{2}^{\ell};\boldsymbol{\xi}_{1}\right]$
and $\boldsymbol{y}_{2}^{\ell}\equiv\left[\boldsymbol{u}_{2}^{\ell};\boldsymbol{v}_{1}^{\ell+1};\boldsymbol{\xi}_{2}\right]$
respectively.

\subsection{Order reduction}

As described in \hyperlink{sec31}{\S3.1}, $\forall i\in\left\{ 1,2\right\} $ and iteration $\ell$, an approximation
of the input data that enters the solver
component $\boldsymbol{m}_{i}$ can be constructed in the reduced dimensional
stochastic space $\Theta_{i}^{\ell}$.

For any square-integrable, vector-valued function $\boldsymbol{u}:\Theta_{i}^{\ell}\rightarrow\mathbb{R}^{n}$,
a reduced gPC approximation $\boldsymbol{u}^{\tilde{p}_{i}}\equiv\boldsymbol{u}^{p,d_{i},\tilde{p}_{i}}$ of order $\tilde{p}_{i}\equiv\tilde{p}_{i}^{\ell}\geq0$
can be formulated as follows. $\forall\boldsymbol{\theta}\in\Theta_{i}^{\ell}$,
\hypertarget{eq322}{}
\begin{equation}
\boldsymbol{u}\left(\boldsymbol{\theta}\right)\approx\boldsymbol{u}^{\tilde{p}_{i}}\left(\boldsymbol{\theta}\right)=\sum_{j=0}^{\tilde{P}_{i}}\tilde{\boldsymbol{u}}_{j}\phi_{i,j}^{\ell}\left(\boldsymbol{\theta}\right)=\tilde{\boldsymbol{U}}\boldsymbol{\phi}_{i}^{\ell,\tilde{p}_{i}}\left(\boldsymbol{\theta}\right),
\end{equation}
where $\left\{ \phi_{i,j}^{\ell}\equiv\phi_{i,j}^{\ell,p}:\Theta_{i}^{\ell}\rightarrow\mathbb{R}\right\} _{j\geq0}$
denotes the set of $\nu_{i}^{\ell}$-orthonormal polynomials,
$\tilde{\boldsymbol{U}}\equiv\tilde{\boldsymbol{U}}^{p,\tilde{p}_{i}}=\left[\begin{array}{ccc}
\tilde{\boldsymbol{u}}_{0} & \cdots & \tilde{\boldsymbol{u}}_{\tilde{P}_{i}}\end{array}\right]\in\mathbb{R}^{n\times\left(\tilde{P}_{i}+1\right)}$ denotes the reduced order gPC coefficient matrix and $\boldsymbol{\phi}_{i}^{\ell,\tilde{p}_{i}}\equiv\boldsymbol{\phi}_{i}^{\ell,p,\tilde{p}_{i}}=\left[\begin{array}{ccc}
\phi_{i,0}^{\ell} & \cdots & \phi_{i,\tilde{P}_{i}}^{\ell}\end{array}\right]^{\mathbf{T}}:\Theta_{i}^{\ell}\rightarrow\mathbb{R}^{\tilde{P}_{i}+1}$ denotes the reduced order basis vector. If the formulation of the gPC approximation $\boldsymbol{u}^{\tilde{p}_{i}}$ is based on an isotropic, total degree truncation of the infinite polynomial series, then $\tilde{p}_{i}$ and $\tilde{P}_{i}$ are related as follows. 
\begin{equation}
\tilde{P}_{i}+1=\frac{\left(\tilde{p}_{i}+d_{i}\right)!}{\tilde{p}_{i}!d_{i}!}.
\end{equation}
Moreover, we assume that $\tilde{p}_{i}\leq p$,
which implies that $\tilde{P}_{i}\leq P$.

Since the coordinate directions in $\Theta_{i}^{\ell}$
are not necessarily statistically independent, we cannot simply compute
the elements of $\boldsymbol{\phi}_{i}^{\ell,\tilde{p}_{i}}$
as products of univariate polynomials. Instead, we propose a
SVD based numerical construction method.

\subsubsection{Reduced order basis construction}

$\forall i\in\left\{ 1,2\right\} $, iteration $\ell$ and degree index $\boldsymbol{\alpha}=\left[\begin{array}{ccc}
\alpha_{1} & \cdots & \alpha_{d_{i}}\end{array}\right]^{\mathbf{T}}\in\mathbb{N}_{0}^{d_{i}}$, let $\tilde{m}_{i,\boldsymbol{\alpha}}^{\ell}:\Theta_{i}^{\ell}\rightarrow\mathbb{R}:\forall\boldsymbol{\theta}=\left[\begin{array}{ccc}
\theta_{1} & \cdots & \theta_{d_{i}}\end{array}\right]\in\Theta_{i}^{\ell}$,
\begin{equation}
\tilde{m}_{i,\boldsymbol{\alpha}}\left(\boldsymbol{\theta}\right)=\prod_{j=1}^{d_{i}}\theta_{j}^{\alpha_{j}}
\end{equation}
denote the monomial function with $\mathrm{deg}\left(\tilde{m}_{i,\boldsymbol{\alpha}}\right)=\left| \boldsymbol{\alpha} \right|=\alpha_{1}+\cdots+\alpha_{d_{i}}$.
The number of such monomial functions with a total degree $\leq\tilde{p}_{i}$
is equal to $\tilde{P}_{i}+1$. Let $\left\{ \boldsymbol{\alpha}_{j}:\left| \boldsymbol{\alpha}_{j}\right| \leq \tilde{p}_{i}\right\}_{j=0}^{\tilde{P}_{i}} $
denote the corresponding set of indices and $\tilde{\boldsymbol{m}}_{i}^{\ell,\tilde{p}_{i}}\equiv\tilde{\boldsymbol{m}}_{i}^{\ell,p,\tilde{p}_{i}}=\left[\begin{array}{ccc}
\tilde{m}_{i,\boldsymbol{\alpha}_{0}} & \cdots & \tilde{m}_{i,\boldsymbol{\alpha}_{\tilde{P}_{i}}}\end{array}\right]:\Theta_{i}^{\ell}\rightarrow\mathbb{R}^{\tilde{P}_{i}+1}$ denote the monomial vector.

The monomial vector $\tilde{\boldsymbol{m}}_{i}^{\ell,\tilde{p}_{i}}$ defines the corresponding Hankel matrix $\boldsymbol{H}_{i}^{\ell,\tilde{p}_{i}}\equiv\boldsymbol{H}_{i}^{\ell,p,\tilde{p}_{i}}:$

\begin{align}
\boldsymbol{H}_{i}^{\ell,\tilde{p}_{i}} & =\int_{\Theta_{i}^{\ell}}\tilde{\boldsymbol{m}}_{i}^{\ell,\tilde{p}_{i}}\left(\boldsymbol{\theta}\right)\tilde{\boldsymbol{m}}_{i}^{\ell,\tilde{p}_{i}}\left(\boldsymbol{\theta}\right)^{\mathbf{T}}\nu_{i}^{\ell}\left(\boldsymbol{\theta}\right)d\boldsymbol{\theta} =\int_{\Xi}\tilde{\boldsymbol{m}}_{i}^{\ell,\tilde{p}_{i}}\left(\boldsymbol{\theta}_{i}^{\ell}\left(\boldsymbol{\xi}\right)\right)\tilde{\boldsymbol{m}}_{i}^{\ell,\tilde{p}_{i}}\left(\boldsymbol{\theta}_{i}^{\ell}\left(\boldsymbol{\xi}\right)\right)^{\mathbf{T}}\mu\left(\boldsymbol{\xi}\right)d\boldsymbol{\xi},
\end{align}
which in turn can be approximated using the global quadrature rule
as follows. $\tilde{\boldsymbol{H}}_{i}^{\ell,\tilde{p}_{i}}\approx\boldsymbol{H}_{i}^{\ell,\tilde{p}_{i}}:$
\begin{align}
\tilde{\boldsymbol{H}}_{i}^{\ell,\tilde{p}_{i}} & =\sum_{j=1}^{Q}w^{\left( j\right)}\tilde{\boldsymbol{m}}_{i}^{\ell,\tilde{p}_{i}}\left(\boldsymbol{\theta}_{i}^{\ell}\left(\boldsymbol{\xi}^{\left( j\right)}\right)\right)\tilde{\boldsymbol{m}}_{i}^{\ell,\tilde{p}_{i}}\left(\boldsymbol{\theta}_{i}^{\ell}\left(\boldsymbol{\xi}^{\left( j\right)}\right)\right)^{\mathbf{T}}.
\end{align}

In general, 
the basis vector $\boldsymbol{\phi}_{i}^{\ell,\tilde{p}_{i}}$
could be computed as follows. $\forall\boldsymbol{\theta}\in \Theta^{\ell}_{i}$, 
\begin{equation}
\boldsymbol{\phi}_{i}^{\ell,\tilde{p}_{i}}\left(\boldsymbol{\theta}\right)=\left( \boldsymbol{L}_{i}^{\ell,\tilde{p}_{i}}\right)^{-1}\tilde{\boldsymbol{m}}_{i}^{\ell,\tilde{p}_{i}}\left(\boldsymbol{\theta}\right),
\end{equation}
where $\boldsymbol{L}_{i}^{\ell,\tilde{p}_{i}}\in\mathbb{R}^{\tilde{P}_{i}+1}$ is the lower triangular matrix which
defines the Cholesky factorization $\boldsymbol{L}_{i}^{\ell,\tilde{p}_{i}}\left(\boldsymbol{L}_{i}^{\ell,\tilde{p}_{i}}\right)^{\mathbf{T}}=\tilde{\boldsymbol{H}}_{i}^{\ell,\tilde{p}_{i}}$. However, the possibility of negative weights in the
quadrature rule can lead to zero or negative eigenvalues in the approximation $\tilde{\boldsymbol{H}}_{i}^{\ell,\tilde{p}_{i}}$. Therefore,
we instead compute the rank-reduced SVD of $\tilde{\boldsymbol{H}}_{i}^{\ell,\tilde{p}_{i}}:$
\begin{equation}
\tilde{\boldsymbol{H}}_{i}^{\ell,\tilde{p}_{i}}=\tilde{\boldsymbol{V}}_{i}^{\ell,\tilde{p}_{i}}\tilde{\boldsymbol{S}_{i}}^{\ell,\tilde{p}_{i}}\tilde{\boldsymbol{\Sigma}}_{i}^{\ell,\tilde{p}_{i}}\left(\tilde{\boldsymbol{V}}_{i}^{\ell,\tilde{p}_{i}}\right)^{\mathbf{T}},
\end{equation}
where $\tilde{\boldsymbol{V}}_{i}^{\ell,\tilde{p}_{i}},\tilde{\boldsymbol{\Sigma}}_{i}^{\ell,\tilde{p}_{i}}$ denote
the usual decomposition matrices and $\tilde{\boldsymbol{S}_{i}}^{\ell,\tilde{p}_{i}}$
denotes a diagonal matrix with $\pm1$ as its diagonal elements. Since
$\tilde{\boldsymbol{H}}_{i}^{\ell,\tilde{p}_{i}}$ is
symmetric, such a decomposition will always exist [\hyperlink{ref29}{29}]. 

Subsequently, the basis vector $\boldsymbol{\phi}_{i}^{\ell,\tilde{p}_{i}}$ is computed as follows.
$\forall\boldsymbol{\theta}\in\tilde{\Theta}_{i}^{\ell}$,
\begin{equation}
\boldsymbol{\phi}_{i}^{\ell,\tilde{p}_{i}}\left(\boldsymbol{\theta}\right)=\left(\tilde{\boldsymbol{\Sigma}}_{i}^{\ell,\tilde{p}_{i}}\right)^{-\frac{1}{2}}\left(\tilde{\boldsymbol{V}}_{i}^{\ell,\tilde{p}_{i}}\right)^{\mathbf{T}}\tilde{\boldsymbol{m}}_{i}^{\ell,\tilde{p}_{i}}\left(\boldsymbol{\theta}\right).
\end{equation}
Here, $\boldsymbol{\phi}_{i}^{\ell,\tilde{p}_{i}}$ satisfies the discrete orthogonality
condition
\begin{equation}
\sum_{j=1}^{Q}w^{\left( j\right)}\boldsymbol{\phi}_{i}^{\ell,\tilde{p}_{i}}\left(\boldsymbol{\theta}_{i}^{\ell}\left(\boldsymbol{\xi}^{\left( j\right)}\right)\right)\boldsymbol{\phi}_{i}^{\ell,\tilde{p}_{i}}\left(\boldsymbol{\theta}_{i}^{\ell}\left(\boldsymbol{\xi}^{\left( j\right)}\right)\right)^{\mathbf{T}}=\tilde{\boldsymbol{S}_{i}}^{\ell,\tilde{p}_{i}},
\end{equation}
and therefore, defines the reduced spectral projection formula
\hypertarget{eq331}{}
\begin{align}
\tilde{\boldsymbol{U}} =\int_{\Theta_{i}^{\ell}}\boldsymbol{u}\left(\boldsymbol{\theta}\right)\boldsymbol{\phi}_{i}^{\ell,\tilde{p}_{i}}\left(\boldsymbol{\theta}\right)^{\mathbf{T}}\tilde{\boldsymbol{S}_{i}}^{\ell,\tilde{p}_{i}}\nu_{i}^{\ell}\left(\boldsymbol{\theta}\right)d\boldsymbol{\theta} & =\int_{\Xi}\boldsymbol{u}\left(\boldsymbol{\theta}_{i}^{\ell}\left(\boldsymbol{\xi}\right)\right)\boldsymbol{\phi}_{i}^{\ell,\tilde{p}_{i}}\left(\boldsymbol{\theta}_{i}^{\ell}\left(\boldsymbol{\xi}\right)\right)^{\mathbf{T}}\tilde{\boldsymbol{S}_{i}}^{\ell,\tilde{p}_{i}}\mu\left(\boldsymbol{\xi}\right)d\boldsymbol{\xi}\nonumber \\
 & \approx\sum_{j=1}^{Q}w^{\left( j\right)}\boldsymbol{u}\left(\boldsymbol{\theta}_{i}^{\ell}\left(\boldsymbol{\xi}^{\left( j\right)}\right)\right)\boldsymbol{\phi}_{i}^{\ell,\tilde{p}_{i}}\left(\boldsymbol{\theta}_{i}^{\ell}\left(\boldsymbol{\xi}^{\left( j\right)}\right)\right)^{\mathbf{T}}\tilde{\boldsymbol{S}_{i}}^{\ell,\tilde{p}_{i}}.
\end{align}

\hyperlink{eq331}{Eq. 3.31} implies that $\left\{ \left(\boldsymbol{\theta}_{i}^{\left( j\right)}=\boldsymbol{\theta}_{i}^{\ell}\left(\boldsymbol{\xi}^{\left( j\right)}\right),w^{\left( j\right)}\right)\right\} _{j=1}^{Q}$
is a quadrature rule of level $\geq\tilde{p}_{i}$ in $\Theta_{i}^{\ell}$.
However, $Q$ depends exponentially on the 
global dimension $s$ and gPC order $p$, implying that the quadrature rule is not 
computationally optimal with respect to the reduced dimension $d_{i}$ and gPC order $\tilde{p}_{i}$. Therefore,
we propose a QR factorization based method to construct the computationally optimal
quadrature rule in $\Theta_{i}^{\ell}$, the size of which depends exponentially on 
$d_{i}$ and $\tilde{p}_{i}$.

\subsubsection{Optimal quadrature rule construction}

$\forall i\in\left\{ 1,2\right\} $ and iteration $\ell$, an optimally sparse quadrature rule 
of level $\tilde{p}_{i}$ in $\Theta_{i}^{\ell}$, denoted as $\left\{ \left(\boldsymbol{\theta}_{i}^{\left( j\right)},\tilde{w}_{i}^{\left( j\right)}\right)\right\} _{j=1}^{Q}$, would
contain the minimum possible number of non-zero weights and be able to numerically integrate
all polynomials with total degree $\leq2\tilde{p}_{i}$, up to machine precision. This requirement is dictated by the 
reduced spectral projection formula in \hyperlink{eq331}{Eq. 3.31}. 

Therefore, if $\boldsymbol{w}=\left[\begin{array}{ccc}
w^{\left(1\right)} & \cdots & w^{\left(Q\right)}\end{array}\right]^{\mathbf{T}}\in\mathbb{R}^{Q}$ and $\tilde{\boldsymbol{w}}_{i}=\left[\begin{array}{ccc}
\tilde{w}_{i}^{\left(1\right)} & \cdots & \tilde{w}^{\left(Q\right)}\end{array}\right]^{\mathbf{T}}\in\mathbb{R}^{Q}$ denote the dense and optimally sparse weight vectors respectively, then $\tilde{\boldsymbol{w}}$
solves the $\ell_{0}-$minimization problem: 
\begin{equation}
\tilde{\boldsymbol{w}}_{i} =\arg\min_{\boldsymbol{\omega}\in\mathbb{R}^{Q}}\left\Vert \boldsymbol{\omega}\right\Vert {}_{0}: \tilde{\boldsymbol{M}}_{i}^{\ell,2\tilde{p}_{i}}\boldsymbol{\omega}=\tilde{\boldsymbol{M}}_{i}^{\ell,2\tilde{p}_{i}}\boldsymbol{w},
\end{equation}
where
\begin{equation}
\tilde{\boldsymbol{M}}_{i}^{\ell,2\tilde{p}_{i}}\equiv\tilde{\boldsymbol{M}}_{i}^{\ell,p,2\tilde{p}_{i}}=\left[\begin{array}{ccc}
\tilde{\boldsymbol{m}}_{i}^{\ell,2\tilde{p}_{i}}\left(\boldsymbol{\theta}_{i}^{\left(1\right)}\right) & \cdots & \tilde{\boldsymbol{m}}_{i}^{\ell,2\tilde{p}_{i}}\left(\boldsymbol{\theta}_{i}^{\left(Q\right)}\right)\end{array}\right]\in\mathbb{R}^{\left(\tilde{N}_{i}+1\right)\times Q}
\end{equation}
denotes the corresponding Vandermonde matrix with $\tilde{N}_{i}+1={\displaystyle \frac{\left(2\tilde{p}_{i}+d_{i}\right)!}{\left(2\tilde{p}_{i}\right)!d_{i}!}}$
rows. In general, $\tilde{\boldsymbol{M}}_{i}^{\ell,2\tilde{p}_{i}}$
may not be full rank and therefore, we can setup a numerically stable
and equivalent $ \ell_{0}-$minimization problem as follows.
\begin{equation}
 \tilde{\boldsymbol{w}}_{i}=\arg\min_{\boldsymbol{\omega}\in\mathbb{R}^{Q}}\left\Vert \boldsymbol{\omega}\right\Vert {}_{0}:
 \left(\boldsymbol{Q}_{r_{i}}^{\ell,2\tilde{p}_{i}}\right)^{\mathbf{T}}\boldsymbol{\omega}=\left(\boldsymbol{Q}_{r_{i}}^{\ell,2\tilde{p}_{i}}\right)^{\mathbf{T}}\boldsymbol{w},
\end{equation}
where $r_{i}$ is the rank of $\tilde{\boldsymbol{M}}_{i}^{\ell,2\tilde{p}_{i}}$
and $\boldsymbol{Q}_{r_{i}}^{\ell,2\tilde{p}_{i}}$ defines the pivoted-QR factorization of
$\left(\tilde{\boldsymbol{M}}_{i}^{\ell,2\tilde{p}_{i}}\right)^{\mathbf{T}}:$
\begin{equation}
\left(\tilde{\boldsymbol{M}}_{i}^{\ell,2\tilde{p}_{i}}\right)^{\mathbf{T}}\boldsymbol{\Pi}^{\ell,2\tilde{p}_{i}}=\boldsymbol{Q}_{r_{i}}^{\ell,2\tilde{p}_{i}}\boldsymbol{R}^{\ell,2\tilde{p}_{i}}.
\end{equation}

Instead of solving the NP-hard $ \ell_{0}-$minimization problem,
we employ a direct approach to construct a 'weakly' optimal quadrature rule, which has at most $r_{i}$ non-zero weights. Therefore, we compute the pivoted-QR
factorization of $\left(\boldsymbol{Q}_{r_{i}}^{\ell,2\tilde{p}_{i}}\right)^{\mathbf{T}}$:
\begin{equation}
\left(\boldsymbol{Q}_{r_{i}}^{\ell,2\tilde{p}_{i}}\right)^{\mathbf{T}}\tilde{\boldsymbol{\Pi}}^{\ell,2\tilde{p}_{i}}=\tilde{\boldsymbol{Q}}^{\ell,2\tilde{p}_{i}}\tilde{\boldsymbol{R}}^{\ell,2\tilde{p}_{i}},
\end{equation}
construct the upper triangular square matrix $\tilde{\boldsymbol{R}}_{r_{i}}^{\ell,2\tilde{p}_{i}}$
using first $r_{i}$ columns of $\tilde{\boldsymbol{R}}^{\ell,2\tilde{p}_{i}}$, and compute the
sparse weight vector as follows. 
\begin{equation}
\tilde{\boldsymbol{w}}_{i}=\tilde{\boldsymbol{\Pi}}^{\ell,2\tilde{p}_{i}}\left[\begin{array}{c}
\left(\tilde{\boldsymbol{R}}_{r_{i}}^{\ell,2\tilde{p}_{i}}\right)^{-1}\tilde{\boldsymbol{R}}^{\ell,2\tilde{p}_{i}}\left(\tilde{\boldsymbol{\Pi}}^{\ell,2\tilde{p}_{i}}\right)^{\mathbf{T}}\boldsymbol{w}\\
\boldsymbol{0}
\end{array}\right].
\end{equation}

If $\tilde{\mathcal{Z}}_{i}\equiv\left\{ 1\leq j\leq Q:\tilde{w}_{i}^{\left( j\right)}\neq0\right\} $
denotes the index set corresponding to the non-zero elements in $\tilde{\boldsymbol{w}}$,
then the reduced spectral projection formula in \hyperlink{eq331}{Eq. 3.31} can be efficiently computed as follows.
\hypertarget{eq338}{}
\begin{equation}
\tilde{\boldsymbol{U}}\approx\sum_{ j\in\tilde{\mathcal{Z}}_{i}}\tilde{w}_{i}^{\left( j\right)}\boldsymbol{u}\left(\boldsymbol{\theta}_{i}^{\left( j\right)}\right)\boldsymbol{\phi}_{i}^{\ell,\tilde{p}_{i}}\left(\boldsymbol{\theta}_{i}^{\left( j\right)}\right)^{\mathbf{T}}\tilde{\boldsymbol{S}_{i}}^{\ell,\tilde{p}_{i}}.
\end{equation}

Using \hyperlink{eq338}{Eq. 3.38}, the global gPC coefficient matrix of $\boldsymbol{u}$
can be subsequently approximated as follows.
\hypertarget{eq339}{}
\begin{align}
\hat{\boldsymbol{U}}^{p}  \approx\hat{\boldsymbol{U}}^{d_{i},\tilde{p}_{i}}\equiv\hat{\boldsymbol{U}}^{p,d_{i},\tilde{p}_{i}} &=\sum_{j=1}^{Q}w^{\left( j\right)}\tilde{\boldsymbol{U}}\boldsymbol{\phi}_{i}^{\ell,\tilde{p}_{i}}\left(\boldsymbol{\theta}_{i}^{\left( j\right)}\right)\boldsymbol{\psi}\left(\boldsymbol{\xi}^{\left( j\right)}\right)^{\mathbf{T}}.
\end{align}

In general, the level and characteristics 
of the global quadrature rule used in approximating the global spectral
projection formula (\hyperlink{eq339}{Eq. 3.39}) can differ from the level and characteristics of the global quadrature rule used in constructing the optimally sparse
quadrature rule. Moreover, a strictly positive weight vector with $r_{i}\leq\tilde{N}_{i}+1$ non-zeros
can be computed using the Nelder-Mead method [\hyperlink{ref30}{30}], in which case,
the upper bound of $\tilde{N}_{i}+1$ on the sparsity of the optimal
weight vector $\tilde{\boldsymbol{w}}_i$ would coincide with the upper bound proved by Tchakaloff
theorem [\hyperlink{ref31}{31}]. Furthermore, the sparsity can also be explicitly controlled
by prescribing a threshold on the diagonal elements of $\boldsymbol{R}^{\ell,2\tilde{p}_{i}}$.

\subsubsection{Selecting the reduced representation}

$\forall i\in\left\{ 1,2\right\} $, we prescribe a
tolerance $\epsilon_{i,\mathrm{ord}}>0$ such that, at each iteration $\ell$, the reduced
order $\tilde{p}_{i}$ is selected as the smallest $k\in\mathbb{N}$,
which satisfies
\hypertarget{eq340}{}
\begin{equation}
\left\Vert \hat{\boldsymbol{U}}_{i}^{\ell,p,d_{i},k+1}-\hat{\boldsymbol{U}}_{i}^{\ell,p,d_{i},k}\right\Vert _{\boldsymbol{\boldsymbol{G}}_{i}}\leq\epsilon_{i,\mathrm{ord}}\left\Vert \hat{\boldsymbol{U}}_{i}^{\ell,p,d_{i},k+1}\right\Vert _{\boldsymbol{\boldsymbol{G}_{i}}},
\end{equation}
where $\left\Vert \cdot\right\Vert _{\boldsymbol{\boldsymbol{G}}_{i}}$ denotes the $\boldsymbol{\boldsymbol{G}}_{i}$-weighted Frobenius
norm, such that $\forall\hat{\boldsymbol{U}}\in\mathbb{R}^{n\times\left(P+1\right)},\left\Vert \hat{\boldsymbol{U}}\right\Vert _{\boldsymbol{\boldsymbol{G}}_{i}}=\sqrt{\mathrm{trace}\left(\hat{\boldsymbol{U}}^{\mathbf{T}}\boldsymbol{\boldsymbol{G}}_{i}\hat{\boldsymbol{U}}\right)}$. 

We propose the following heuristic to select $\tilde{p}_{i}$. $\tilde{p}_{i}$ is initialized to $0$ at $\ell=0$, and  
subsequently incremented by $1$ at any $\ell>0$, if $\left\Vert \hat{\boldsymbol{U}}_{i}^{\ell,p,d_{i},\tilde{p}_{i}+1}-\hat{\boldsymbol{U}}_{i}^{\ell,p,d_{i},\tilde{p}_{i}}\right\Vert _{\boldsymbol{\boldsymbol{G}}_{i}}>\epsilon_{i,\mathrm{ord}}\left\Vert \hat{\boldsymbol{U}}_{i}^{\ell,p,d_{i},\tilde{p}_{i}+1}\right\Vert _{\boldsymbol{\boldsymbol{G}_{i}}}$.
Therefore, by choosing an appropriate value for the tolerance $\epsilon_{i,\mathrm{ord}}$, we can guarantee that $\tilde{p}_{i}<p$ and therefore, a reduction of the gPC approximation order.
Moreover, the requirement of computing $\hat{\boldsymbol{U}}_{i}^{\ell,p,d_{i},\tilde{p}_{i}+1}$ implies that
the level of the constructed optimal quadrature rule must be $\geq\tilde{p}_{i}+1$. 

\hyperlink{thm2}{Theorem 2} proves an important relation between the
tolerance $\epsilon_{i,\mathrm{ord}}$ and the upper
bound on the error incurred by the reduced order approximation defined in
\hyperlink{eq322}{Eq. 3.22}.
\subsubsection*{Theorem 2:}
\hypertarget{thm2}{}
\emph{
$\forall i\in\left\{1,2\right\}$ and iteration $\ell$, let $\boldsymbol{u},\tilde{\boldsymbol{u}}^{\tilde{p}_{i}}:\Theta_{i}^{\ell}\rightarrow\mathbb{R}^{n}$
denote an infinitely regular vector valued function in
$\Theta_{i}^{\ell}$ and its reduced gPC approximation
of order $\tilde{p}_{i}\geq0$ respectively. Given a positive-definite Gramian 
$\boldsymbol{G}\in\mathbb{R}^{n\times n}$, if $\boldsymbol{u}^{d_{i}}\equiv\boldsymbol{u}^{p,d_{i}}$
and $\boldsymbol{u}^{d_{i},\tilde{p}_{i}}\equiv\boldsymbol{u}^{p,d_{i},\tilde{p}_{i}}$ denote
the respective global gPC approximations of the composite functions
$\boldsymbol{u}\circ\boldsymbol{\theta}_{i}^{\ell}$ and
$\tilde{\boldsymbol{u}}^{\tilde{p}_{i}}\circ\boldsymbol{\theta}_{i}^{\ell}$,
then $\exists\chi,\chi^{*}>0,\rho^{*}>1$, such that $\forall p,\tilde{p}_{i}\geq 0$,
\hypertarget{eq341}{}
\begin{align}
\sqrt{\int_{\Xi}\left\Vert \boldsymbol{u}^{d_{i}}\left(\boldsymbol{\xi}\right)-\boldsymbol{u}^{d_{i},\tilde{p}_{i}}\left(\boldsymbol{\xi}\right)\right\Vert _{\boldsymbol{\boldsymbol{G}}}^{2}\mu\left(\boldsymbol{\xi}\right)d\boldsymbol{\xi}} & \leq\chi\epsilon_{i,\mathrm{ord}}\sqrt{\int_{\Xi}\left\Vert \boldsymbol{u}^{d_{i}}\left(\boldsymbol{\xi}\right)\right\Vert _{\boldsymbol{\boldsymbol{G}}}^{2}\mu\left(\boldsymbol{\xi}\right)d\boldsymbol{\xi}} +\chi^{*}\rho^{-p}
\end{align}
for some $\rho\geq\rho^{*}$. }
\subsubsection*{Proof:}
Since $\boldsymbol{u}$ is infinitely regular in $\Theta_{i}^{\ell}$,
$\boldsymbol{u}\circ\boldsymbol{\theta}_{i}^{\ell}$ is
infinitely regular in $\Xi$. Therefore, from the Cameron-Martin
theorem, $\exists\chi^{*},\tilde{\chi}>0,\rho^{*}>1$, such that for any positive-definite
$\boldsymbol{G}\in\mathbb{R}^{n\times n}$, the approximation error
between $\boldsymbol{u}^{d_{i}}$ and $\boldsymbol{u}^{d_{i},\tilde{p}_{i}}$
has the following upper bound.

\hypertarget{eq342}{}
\begin{align}
& \sqrt{\int_{\Xi}\left\Vert \boldsymbol{u}^{d_{i}}\left(\boldsymbol{\xi}\right)-\tilde{\boldsymbol{u}}^{d_{i},\tilde{p}_{i}}\left(\boldsymbol{\xi}\right)\right\Vert _{\boldsymbol{\boldsymbol{G}}}^{2}\mu\left(\boldsymbol{\xi}\right)d\boldsymbol{\xi}}\nonumber \\
& \leq\tilde{\chi}\sqrt{\int_{\Xi}\left\Vert \tilde{\boldsymbol{u}}^{d_{i},\tilde{p}_{i}+1}\left(\boldsymbol{\xi}\right)-\tilde{\boldsymbol{u}}^{d_{i},\tilde{p}_{i}}\left(\boldsymbol{\xi}\right)\right\Vert _{\boldsymbol{\boldsymbol{G}}}^{2}\mu\left(\boldsymbol{\xi}\right)d\boldsymbol{\xi}}+\chi^{*}\rho^{-p}.\nonumber \\
& = \tilde{\chi}\sqrt{\int_{\Xi}\left\Vert \left(\hat{\boldsymbol{U}}^{d_{i},\tilde{p}_{i}+1}-\hat{\boldsymbol{U}}^{d_{i},\tilde{p}_{i}}\right)\boldsymbol{\psi}\left(\boldsymbol{\xi}\right)\right\Vert _{\boldsymbol{\boldsymbol{G}}}^{2}\mu\left(\boldsymbol{\xi}\right)d\boldsymbol{\xi}}+\chi^{*}\rho^{-p}\nonumber\\
& = \tilde{\chi}\sqrt{\int_{\Xi}\left\Vert \boldsymbol{G}^{\frac{1}{2}}\left(\hat{\boldsymbol{U}}^{d_{i},\tilde{p}_{i}+1}-\hat{\boldsymbol{U}}^{d_{i},\tilde{p}_{i}}\right)\boldsymbol{\psi}\left(\boldsymbol{\xi}\right)\right\Vert _{2}^{2}\mu\left(\boldsymbol{\xi}\right)d\boldsymbol{\xi}}+\chi^{*}\rho^{-p} \nonumber \\
& \leq \tilde{\chi}\left\Vert \boldsymbol{G}^{\frac{1}{2}}\left(\hat{\boldsymbol{U}}^{d_{i},\tilde{p}_{i}+1}-\hat{\boldsymbol{U}}^{d_{i},\tilde{p}_{i}}\right)\right\Vert _{2}\sqrt{\int_{\Xi}\left\Vert \boldsymbol{\psi}\left(\boldsymbol{\xi}\right)\right\Vert _{2}^{2}\mu\left(\boldsymbol{\xi}\right)d\boldsymbol{\xi}}+\chi^{*}\rho^{-p} \nonumber \\
& \leq \chi\left\Vert \hat{\boldsymbol{U}}^{d_{i},\tilde{p}_{i}+1}-\hat{\boldsymbol{U}}^{d_{i},\tilde{p}_{i}}\right\Vert _{\boldsymbol{\boldsymbol{G}}}+\chi^{*}\rho^{-p}
\end{align}
for some $\chi>0,\rho\geq\rho^{*}$. By substituting \hyperlink{eq340}{Eq. 3.40} into \hyperlink{eq342}{Eq. 3.42}, we arrive at \hyperlink{eq341}{Eq. 3.41} $\square$

\subsection{Algorithm and computational cost}

\begin{algorithm}
\hypertarget{alg2}{}
\caption{Reduced NISP based uncertainty propagation for a two-module multi-physics system}
\SetKwInOut{Input}{inputs}\SetKwInOut{Output}{outputs}
\SetKwRepeat{Repeat}{repeat}{}
\DontPrintSemicolon
\Input{$\mu_{1}$, $\mu_{2}$, $p\geq0$, $q\geq p$, $\epsilon_{1,\mathrm{dim}}$,
$\epsilon_{2,\mathrm{dim}}$, $\epsilon_{1,\mathrm{ord}}$,
$\epsilon_{2,\mathrm{ord}}$, $\hat{\boldsymbol{U}}_{1}^{0}$, $\hat{\boldsymbol{U}}_{2}^{0}$} 
\Output{$\hat{\boldsymbol{U}}_{1}$, $\hat{\boldsymbol{U}}_{2}$}
\textbf{precompute}: $\left\{ \left(w^{\left( j\right)},\boldsymbol{\psi}\left(\boldsymbol{\xi}^{\left( j\right)}\right)\right)\right\} _{j=1}^{Q}$,
$\hat{\boldsymbol{\Xi}}_{1}$, $\hat{\boldsymbol{\Xi}}_{2}$

\textbf{$\ell\leftarrow0$}, \textbf{$\tilde{p}_{1}\leftarrow0$}, \textbf{$\tilde{p}_{2}\leftarrow0$}

\Repeat{(contd.)}{

\SetKwBlock{Begin}{dimension reduction}{end}
\Begin{
\Input{$\hat{\boldsymbol{U}}_{1}^{\ell},\hat{\boldsymbol{U}}_{2}^{\ell},\hat{\boldsymbol{\Xi}}_{1},\epsilon_{1,\mathrm{dim}}$}

\Output{$\bar{\boldsymbol{u}}_{1,1}^{\ell},\tilde{\boldsymbol{U}}_{1,1}^{\ell},\bar{\boldsymbol{u}}_{2,1}^{\ell},\tilde{\boldsymbol{U}}_{2,1}^{\ell},\bar{\boldsymbol{\xi}}_{1,1}^{\ell},\tilde{\boldsymbol{\Xi}}_{1,1}^{\ell},\left\{ \boldsymbol{\theta}_{1}^{\left( j\right)}\right\} _{j=1}^{Q}$}
}
\SetKwBlock{Begin}{reduced basis/quadrature construction}{end}
\Begin{
\Input{$\left\{ \left(\boldsymbol{\theta}_{1}^{\left( j\right)},w^{\left( j\right)}\right)\right\} _{j=1}^{Q}$
$,\tilde{p}_{1}$}

\Output{$\left\{ \left(\boldsymbol{\phi}_{1}^{\ell,\tilde{p}_{1}}\left(\boldsymbol{\theta}_{1}^{\left( j\right)}\right),\boldsymbol{\phi}_{1}^{\ell,\tilde{p}_{1}+1}\left(\boldsymbol{\theta}_{1}^{\left( j\right)}\right),\tilde{w}_{1}^{\left( j\right)}\right)\right\} _{j=1}^{Q}$, 
$\boldsymbol{S}_{1}^{\ell,\tilde{p}_{1}}$,
$\boldsymbol{S}_{1}^{\ell,\tilde{p}_{1}+1}$, $\tilde{\mathcal{Z}}_{1}$}
}

\textbf{}$\tilde{\boldsymbol{U}}_{1}^{\tilde{p}_{1}}\leftarrow\boldsymbol{0}$,
$\tilde{\boldsymbol{U}}_{1}^{\tilde{p}_{1}+1}\leftarrow\boldsymbol{0}$

\For{$ j\in\tilde{\mathcal{Z}}_{1}$}{

\textbf{$\boldsymbol{u}_{1}\leftarrow\boldsymbol{m}_{1}\left(\bar{\boldsymbol{u}}_{1,1}^{\ell}+\tilde{\boldsymbol{U}}_{1,1}^{\ell}\boldsymbol{\theta}_{1}^{\left( j\right)},\bar{\boldsymbol{u}}_{2,1}^{\ell}+\tilde{\boldsymbol{U}}_{2,1}^{\ell}\boldsymbol{\theta}_{1}^{\left( j\right)},\bar{\boldsymbol{\xi}}_{1,1}^{\ell}+\tilde{\boldsymbol{\Xi}}_{1,1}^{\ell}\boldsymbol{\theta}_{1}^{\left( j\right)}\right)$}

\textbf{}$\tilde{\boldsymbol{U}}_{1}^{\tilde{p}_{1}}\leftarrow\tilde{\boldsymbol{U}}_{1}^{\tilde{p}_{1}}+\tilde{w}_{1}^{\left( j\right)}\boldsymbol{u}_{1}\boldsymbol{\phi}_{1}^{\ell,\tilde{p}_{1}}\left(\boldsymbol{\theta}_{1}^{\left( j\right)}\right)^{\mathbf{T}}\boldsymbol{S}_{1}^{\ell,\tilde{p}_{1}}$

\textbf{$\tilde{\boldsymbol{U}}_{1}^{\tilde{p}_{1}+1}\leftarrow\tilde{\boldsymbol{U}}_{1}^{\tilde{p}_{1}+1}+\tilde{w}_{1}^{\left( j\right)}\boldsymbol{u}_{1}\boldsymbol{\phi}_{1}^{\ell,\tilde{p}_{1}+1}\left(\boldsymbol{\theta}_{1}^{\left( j\right)}\right)^{\mathbf{T}}\boldsymbol{S}_{1}^{\ell,\tilde{p}_{1}+1}$}

}

\textbf{}$\hat{\boldsymbol{U}}_{1}^{\tilde{p}_{1}}\leftarrow\boldsymbol{0},\hat{\boldsymbol{U}}_{1}^{\tilde{p}_{1}+1}\leftarrow\boldsymbol{0}$

\For{$ j\leftarrow 1$ \KwTo $Q$}{

\textbf{$\hat{\boldsymbol{U}}_{1}^{\tilde{p}_{1}}\leftarrow\hat{\boldsymbol{U}}_{1}^{\tilde{p}_{1}}+w^{\left( j\right)}\tilde{\boldsymbol{U}}_{1}^{\tilde{p}_{1}}\boldsymbol{\phi}_{1}^{\ell,\tilde{p}_{1}}\left(\boldsymbol{\theta}_{1}^{\left( j\right)}\right)\boldsymbol{\psi}\left(\boldsymbol{\xi}^{\left( j\right)}\right)^{\mathbf{T}}$}
\textbf{$\hat{\boldsymbol{U}}_{1}^{\tilde{p}_{1}+1}\leftarrow\hat{\boldsymbol{U}}_{1}^{\tilde{p}_{1}+1}+w^{\left( j\right)}\tilde{\boldsymbol{U}}_{1}^{\tilde{p}_{1}+1}\boldsymbol{\phi}_{1}^{\ell,\tilde{p}_{1}+1}\left(\boldsymbol{\theta}_{1}^{\left( j\right)}\right)\boldsymbol{\psi}\left(\boldsymbol{\xi}^{\left( j\right)}\right)^{\mathbf{T}}$}
}
\If{$\left\Vert \hat{\boldsymbol{U}}_{1}^{\tilde{p}_{1}+1}-\hat{\boldsymbol{U}}_{1}^{\tilde{p}_{1}}\right\Vert _{\boldsymbol{\boldsymbol{G}}_{1}}>\epsilon_{1,\mathrm{ord}}\left\Vert \hat{\boldsymbol{U}}_{1}^{\tilde{p}_{1}+1}\right\Vert _{\boldsymbol{\boldsymbol{G}}_{1}}$}
{$\tilde{p}_{1}\leftarrow\tilde{p}_{1}+1$}

$\hat{\boldsymbol{U}}_{1}^{\ell+1}\leftarrow\hat{\boldsymbol{U}}_{1}^{\tilde{p}_{1}}$

}
\end{algorithm}
\RestyleAlgo{boxed}
\begin{algorithm}
\SetKwInOut{Input}{inputs}\SetKwInOut{Output}{outputs}
\SetKwRepeat{Repeat}{}{until}
\DontPrintSemicolon
 \emph{(contd.)}\Repeat{$\hat{\boldsymbol{U}}_{1}^{\ell}$, $\hat{\boldsymbol{U}}_{2}^{\ell}$ $\mathrm{converge}$}{

\SetKwBlock{Begin}{dimension reduction}{end}
\Begin{

\Input{$\hat{\boldsymbol{U}}_{2}^{\ell},\hat{\boldsymbol{U}}_{1}^{\ell+1},\hat{\boldsymbol{\Xi}}_{2},\epsilon_{2,\mathrm{dim}}$}

\Output{$\bar{\boldsymbol{u}}_{1,2}^{\ell},\tilde{\boldsymbol{U}}_{1,2}^{\ell},\bar{\boldsymbol{u}}_{2,2}^{\ell},\tilde{\boldsymbol{U}}_{2,2}^{\ell},\bar{\boldsymbol{\xi}}_{2,2}^{\ell},\tilde{\boldsymbol{\Xi}}_{2,2}^{\ell},\left\{ \boldsymbol{\theta}_{2}^{\left( j\right)}\right\} _{j=1}^{Q}$}
}

\SetKwBlock{Begin}{reduced basis/quadrature construction}{end}
\Begin{

\Input{$\left\{ \left(\boldsymbol{\theta}_{2}^{\left( j\right)},w^{\left( j\right)}\right)\right\} _{j=1}^{Q}$
$,\tilde{p}_{2}$}

\Output{$\left\{ \left(\boldsymbol{\phi}_{2}^{\ell,\tilde{p}_{2}}\left(\boldsymbol{\theta}_{2}^{\left( j\right)}\right),\boldsymbol{\phi}_{2}^{\ell,\tilde{p}_{2}+1}\left(\boldsymbol{\theta}_{2}^{\left( j\right)}\right),\tilde{w}_{2}^{\left( j\right)}\right)\right\} _{j=1}^{Q}$, 
$\boldsymbol{S}_{2}^{\ell,\tilde{p}_{2}}$,
$\boldsymbol{S}_{2}^{\ell,\tilde{p}_{2}+1}$, $\tilde{\mathcal{Z}}_{2}$}
}
\textbf{}$\tilde{\boldsymbol{U}}_{2}^{\tilde{p}_{2}}\leftarrow\boldsymbol{0}$,
$\tilde{\boldsymbol{U}}_{2}^{\tilde{p}_{2}+1}\leftarrow\boldsymbol{0}$

\For{$ j\in\tilde{\mathcal{Z}}_{2}$}{

\textbf{$\boldsymbol{u}_{2}\leftarrow\boldsymbol{m}_{2}\left(\bar{\boldsymbol{u}}_{2,2}^{\ell}+\tilde{\boldsymbol{U}}_{2,2}^{\ell}\boldsymbol{\theta}_{2}^{\left( j\right)},\bar{\boldsymbol{u}}_{1,2}^{\ell}+\tilde{\boldsymbol{U}}_{1,2}^{\ell}\boldsymbol{\theta}_{2}^{\left( j\right)},\bar{\boldsymbol{\xi}}_{2,2}^{\ell}+\tilde{\boldsymbol{\Xi}}_{2,2}^{\ell}\boldsymbol{\theta}_{2}^{\left( j\right)}\right)$}

\textbf{}$\tilde{\boldsymbol{U}}_{2}^{\tilde{p}_{2}}\leftarrow\tilde{\boldsymbol{U}}_{2}^{\tilde{p}_{2}}+\tilde{w}_{2}^{\left( j\right)}\boldsymbol{u}_{2}\boldsymbol{\phi}_{2}^{\ell,\tilde{p}_{2}}\left(\boldsymbol{\theta}_{2}^{\left( j\right)}\right)^{\mathbf{T}}\boldsymbol{S}_{2}^{\ell,\tilde{p}_{2}}$

\textbf{$\tilde{\boldsymbol{U}}_{2}^{\tilde{p}_{2}+1}\leftarrow\tilde{\boldsymbol{U}}_{2}^{\tilde{p}_{2}+1}+\tilde{w}_{2}^{\left( j\right)}\boldsymbol{u}_{2}\boldsymbol{\phi}_{2}^{\ell,\tilde{p}_{2}+1}\left(\boldsymbol{\theta}_{2}^{\left( j\right)}\right)^{\mathbf{T}}\boldsymbol{S}_{2}^{\ell,\tilde{p}_{2}+1}$}

}

\textbf{}$\hat{\boldsymbol{U}}_{2}^{\tilde{p}_{2}}\leftarrow\boldsymbol{0},\hat{\boldsymbol{U}}_{2}^{\tilde{p}_{2}+1}\leftarrow\boldsymbol{0}$

\For {$ j\leftarrow 1$ \KwTo $Q$}{

\textbf{$\hat{\boldsymbol{U}}_{2}^{\tilde{p}_{2}}\leftarrow\hat{\boldsymbol{U}}_{2}^{\tilde{p}_{2}}+w^{\left( j\right)}\tilde{\boldsymbol{U}}_{2}^{\tilde{p}_{2}}\boldsymbol{\phi}_{2}^{\ell,\tilde{p}_{2}}\left(\boldsymbol{\theta}_{2}^{\left( j\right)}\right)\boldsymbol{\psi}\left(\boldsymbol{\xi}^{\left( j\right)}\right)^{\mathbf{T}}$}

\textbf{$\hat{\boldsymbol{U}}_{2}^{\tilde{p}_{2}+1}\leftarrow\hat{\boldsymbol{U}}_{2}^{\tilde{p}_{2}+1}+w^{\left( j\right)}\tilde{\boldsymbol{U}}_{2}^{\tilde{p}_{2}+1}\boldsymbol{\phi}_{2}^{\ell,\tilde{p}_{2}+1}\left(\boldsymbol{\theta}_{2}^{\left( j\right)}\right)\boldsymbol{\psi}\left(\boldsymbol{\xi}^{\left( j\right)}\right)^{\mathbf{T}}$}

}

\If{$\left\Vert \hat{\boldsymbol{U}}_{2}^{\tilde{p}_{2}+1}-\hat{\boldsymbol{U}}_{2}^{\tilde{p}_{2}}\right\Vert _{\boldsymbol{\boldsymbol{G}}_{2}}>\epsilon_{2,\mathrm{ord}}\left\Vert \hat{\boldsymbol{U}}_{2}^{\tilde{p}_{2}+1}\right\Vert _{\boldsymbol{\boldsymbol{G}}_{2}}$}
{$\tilde{p}_{2}\leftarrow\tilde{p}_{2}+1$}

$\hat{\boldsymbol{U}}_{2}^{\ell+1}\leftarrow\hat{\boldsymbol{U}}_{2}^{\tilde{p}_{2}}$

\textbf{$\ell\leftarrow \ell+1$}

}
\end{algorithm}

The proposed reduced NISP based uncertainty propagation method is described in \hyperlink{alg2}{Algorithm 2}. Let $\tilde{Q}_{1}$ denote the size of
the optimal quadrature rule corresponding to the reduced dimension $d_{1}$ and reduced order $\tilde{p}_{1}$, and let $\tilde{Q}_{2}$ denote the size of the
optimal quadrature rule corresponding to the 
reduced dimension $d_{2}$ and reduced order $\tilde{p}_{2}$. Therefore, retaining the assumption that the computational costs are dominated by the repeated execution of module operators $\boldsymbol{m}_{1}$
and $\boldsymbol{m}_{2}$, the computational cost of the reduced NISP method
\begin{equation}
C_{r}\approx\mathcal{O}\left(\bar{\mathcal{C}}_{1}\tilde{Q}_{1}+\bar{\mathcal{C}}_{2}\tilde{Q}_{2}\right).
\end{equation}

would grow exponentially with respect to the reduced dimensions $d_{1}.d_{2}$ and orders $\tilde{p}_{1},\tilde{p}_{2}$. Therefore,
the proposed reduced NISP method indeed mitigates the curse of dimensionality associated with the standard NISP method.

\subsection{Error analysis}
\hypertarget{sec34}{}

$\forall i\in\left\{ 1,2\right\} $, let $\varepsilon_{i}:\forall\boldsymbol{\xi}\in\Xi$,
\begin{equation}
\varepsilon_{i}\left(\boldsymbol{\xi}\right)=\left\Vert \boldsymbol{u}_{i}\left(\boldsymbol{\xi}\right)-\boldsymbol{u}_{i}^{\ell,p,d_{i},\tilde{p}_{i}}\left(\boldsymbol{\xi}\right)\right\Vert _{\boldsymbol{G}_{i}}
\end{equation}
denote the mean-square erro between the component solution $\boldsymbol{u}_{i}$ and its corresponding reduced gPC approximation $\boldsymbol{u}_{i}^{\ell,p,d_{i},\tilde{p}_{i}}$.

Subsequently, $\varepsilon_{i}$
can be decompoed as a sum of individual error terms as follows. 
\begin{equation}
\varepsilon_{i}=\varepsilon_{i,\mathrm{BGS}}+\varepsilon_{i,\mathrm{gPC}}+\varepsilon_{i,\mathrm{dim}}+\varepsilon_{i,\mathrm{ord}},
\end{equation}
where $\varepsilon_{i,\mathrm{BGS}}$ denotes the convergence
error, $\varepsilon_{i,\mathrm{gPC}}$ denotes the
gPC truncation error, $\varepsilon_{i,\mathrm{dim}}$
denotes the reduced dimension approximation error and $\varepsilon_{i,\mathrm{ord}}$
denotes the reduced order approximation error. Using the triangle inequality property of norms, an asymptotic upper bound for each  constituent error term
can be formulated as follows. $\forall\boldsymbol{\xi}\in\Xi$,

\hypertarget{eq346}{}
\begin{align}
\varepsilon_{i}\left(\boldsymbol{\xi}\right)= & \left\Vert \boldsymbol{u}_{i}\left(\boldsymbol{\xi}\right)-\boldsymbol{u}_{i}^{\ell}\left(\boldsymbol{\xi}\right)+\boldsymbol{u}_{i}^{\ell}\left(\boldsymbol{\xi}\right)-\boldsymbol{u}_{i}^{\ell,p}\left(\boldsymbol{\xi}\right)+\boldsymbol{u}_{i}^{\ell,p}\left(\boldsymbol{\xi}\right)-\boldsymbol{u}_{i}^{\ell,p,d_{i}}\left(\boldsymbol{\xi}\right)+\boldsymbol{u}_{i}^{\ell,p,d_{i}}\left(\boldsymbol{\xi}\right)\right.\nonumber\\
& \left.-\boldsymbol{u}_{i}^{\ell,p,d_{i},\tilde{p}_{i}}\left(\boldsymbol{\xi}\right)\right\Vert _{\boldsymbol{G}_{i}}\nonumber \\
\leq & \underset{\varepsilon_{i,\mathrm{BGS}}\left(\boldsymbol{\xi}\right)\leq\mathcal{O}\left(\eta{}^{-\ell}\right)}{\underbrace{\left\Vert \boldsymbol{u}_{i}\left(\boldsymbol{\xi}\right)-\boldsymbol{u}_{i}^{\ell}\left(\boldsymbol{\xi}\right)\right\Vert _{\boldsymbol{G}_{i}}}}+\underset{\varepsilon_{i,\mathrm{gPC}}\left(\boldsymbol{\xi}\right)\leq\mathcal{O}\left(\rho^{-p}\right)}{\underbrace{\left\Vert \boldsymbol{u}_{i}^{\ell}\left(\boldsymbol{\xi}\right)-\boldsymbol{u}_{i}^{\ell,p}\left(\boldsymbol{\xi}\right)\right\Vert _{\boldsymbol{G}_{i}}}}+\underset{\varepsilon_{i,\mathrm{dim}}\left(\boldsymbol{\xi}\right)\leq\mathcal{O}\left(\epsilon_{i,\mathrm{dim}}\right)}{\underbrace{\left\Vert \boldsymbol{u}_{i}^{\ell,p}\left(\boldsymbol{\xi}\right)-\boldsymbol{u}_{i}^{\ell,p,d_{i}}\left(\boldsymbol{\xi}\right)\right\Vert _{\boldsymbol{G}_{i}}}}\nonumber \\
 & +\underset{\varepsilon_{i,\mathrm{ord}}\left(\boldsymbol{\xi}\right)\leq\mathcal{O}\left(\epsilon_{i,\mathrm{ord}}\right)+\mathcal{O}\left(\tilde{\rho}^{-p}\right)}{\underbrace{\left\Vert \boldsymbol{u}_{i}^{\ell,p,d_{i}}\left(\boldsymbol{\xi}\right)-\boldsymbol{u}_{i}^{\ell,p,d_{i},\tilde{p}_{i}}\left(\boldsymbol{\xi}\right)\right\Vert _{\boldsymbol{G}_{i}}}}.
\end{align}

In the standard NISP method, the asymptotic upper bound on the approximation error
would simply be the sum of the first two terms on the right hand side
of \hyperlink{eq346}{Eq. 3.46}, implying that $\boldsymbol{u}_{i}^{\ell,p}$
would converge to $\boldsymbol{u}_{i}$ as $\ell,p\rightarrow\infty$.
However, in the reduced NISP method, $\varepsilon_{i}$
would converge to a non-zero quantity, and have an asymptotic upper bound of $\mathcal{O}\left(\epsilon_{i,\mathrm{dim}}\right)+\mathcal{O}\left(\epsilon_{i,\mathrm{ord}}\right)$.

This analysis suggests that the approximation error in $\boldsymbol{u}_{i}^{\ell,p,d_{i},\tilde{p}_{i}}$ can be explicitly controlled by choosing the tolerances $\epsilon_{i,\mathrm{dim}}$ and $\epsilon_{i,\mathrm{ord}}$ appropriately.

\subsection{Selecting the tolerance values}
In practice, since the exact bounds on the approximation error are not known apriori, the tolerance values are selected based on the results computed in a preliminary study, wherein a lower fidelity multi-physics model is used and the various tolerance values can be tested. In our case, lower fidelity translates to a coarser spatial discretization of the original steady state coupled PDE system to formulate the coupled algebraic system. An illustration of this is provided in the numerical experiments in \hyperlink{sec4}{\S4}

\section{Numerical examples}
\hypertarget{sec4}{}
We will now demonstrate and compare the performance of the standard and reduced NISP methods using two numerical examples.

\subsection{Poisson problem}
\hypertarget{sec41}{}
The Poisson problem is a widely used numerical example for benchmarking
uncertainty propagation methods [\hyperlink{ref}{32}, \hyperlink{ref33}{33}]. We consider a system
of Poisson equations across two coupled domains with uncertain diffusion
coefficients.

\subsubsection{Model setup}

Let $\Omega_{1}\equiv\left(-1,0\right)_{x_{1}}\times\left(0,1\right)_{x_{2}}$
and $\Omega_{2}\equiv\left(0,1\right)_{x_{1}}\times\left(0,1\right)_{x_{2}}$
denote the non-overlapping spatial domains , $\Gamma_{12}\equiv\left\{ 0\right\} _{x_{1}}\times\left(0,1\right)_{x_{2}}$
denote the interface and $\boldsymbol{x}=\left[\begin{array}{cc}
x_{1} & x_{2}\end{array}\right]^{\mathbf{T}}$ denote a point in $\Omega_{1}\cup\Omega_{2}\cup\Gamma_{12}$. The
solution fields $u_{1}$ and $u_{2}$ are governed by the following PDE system. $\forall \boldsymbol{\xi} \in \Xi$,
\begin{align}
\boldsymbol{\nabla}^{\mathbf{T}}\left(a_{1}\left(\boldsymbol{x},\boldsymbol{\xi}_{1}\left(\boldsymbol{\xi}\right)\right)\boldsymbol{\nabla}u_{1}\left(\boldsymbol{x},\boldsymbol{\xi}\right)\right)+b_{1} & =0, & \boldsymbol{x}\in\Omega_{1},\nonumber \\
\boldsymbol{\nabla}^{\mathbf{T}}\left(a_{2}\left(\boldsymbol{x},\boldsymbol{\xi}_{2}\left(\boldsymbol{\xi}\right)\right)\boldsymbol{\nabla}u_{2}\left(\boldsymbol{x},\boldsymbol{\xi}\right)\right)+b_{2} & =0, & \boldsymbol{x}\in\Omega_{2},
\end{align}
with the interface condition 
\begin{align}
u_{1}\left(\boldsymbol{x},\boldsymbol{\xi}\right)-u_{2}\left(\boldsymbol{x},\boldsymbol{\xi}\right) & =0, & \boldsymbol{x}\in\Gamma_{12},
\end{align}
and homogenous Dirichlet boundary conditions on all other boundaries.
\hyperlink{fig1}{Figure 1} illustrates the computational domain.
\begin{figure}
    \hypertarget{fig1}{}
    \centering
    \includegraphics[bb=125bp 340bp 475bp 550bp,scale=0.7]{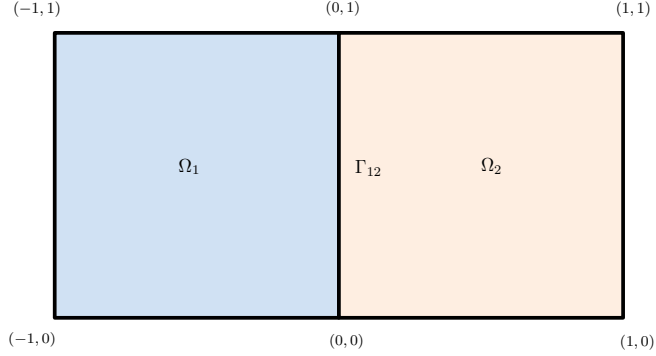}  
    \caption{Computational domain for the Poisson problem.}
    \label{fig:fig1}
\end{figure}

The stochastic diffusion coefficients $a_{1}$ and $a_{2}$ are modeled
using the following KL expansions. $\forall i\in\left\{ 1,2\right\} ,x\in\Omega_{i},\boldsymbol{\xi}_{i}\in\Xi_{i}\equiv\left[-1,1\right]^{s_{i}}$,
\begin{equation}
a_{i}\left(\boldsymbol{x},\boldsymbol{\xi}_{i}\right)=\bar{a}_{i}+\sqrt{3}\delta_{i}\sum_{j=1}^{s_{i}}\gamma_{i,j}\left(\boldsymbol{x}\right)\xi_{ij},
\end{equation}
where $\bar{a}_{i}$ denotes the mean of $a_{i}$ and $\left\{ \xi_{ij}\sim U\left[-1,1\right]\right\}_{j=1}^{s_{i}} $
are $i.i.d.$ random variables. Moreover, we assume that $a_{i}$ has an exponential
covariance kernel
\begin{align}
C_{a_{i}}\left(\boldsymbol{x},\boldsymbol{y}\right)=\delta_{i}^{2}\exp\left(-\frac{\left\Vert \boldsymbol{x}-\boldsymbol{y}\right\Vert _{1}}{l_{i}}\right),& & \boldsymbol{x},\boldsymbol{y}\in\Omega_{i},
\end{align}
where $\delta_{i}$ denotes the coefficient of variation and $l_{i}$
denotes the correlation length of $a_{i}$. The analytic expressions for
$\left\{\gamma_{i,j}\right\}_{j>0}$ are provided in \hyperlink{appA}{Appendix A}. 

The governing PDE system is spatially discretized using bilinear finite
elements [\hyperlink{ref34}{34}] and $m\times m$ equispaced nodes in each subdomain.
We therefore obtain linear equations of the form $\boldsymbol{A}_{1}\left(\boldsymbol{\xi}_{1}\right)\boldsymbol{u}_{1}^{\prime}=\boldsymbol{b}_{1}$
and $\boldsymbol{A}_{2}\left(\boldsymbol{\xi}_{2}\right)\boldsymbol{u}_{2}^{\prime}=\boldsymbol{b}_{2}$
where $\boldsymbol{u}_{1}^{\prime},\boldsymbol{u}_{2}^{\prime}\in\mathbb{R}^{m^{2}}$
are the vectors of nodal values, $\boldsymbol{A}_{1},\boldsymbol{A}_{2}\in\mathbb{R}^{m^{2}\times m^{2}}$
are the stiffness matrices and $\boldsymbol{b}_{1},\boldsymbol{b}_{2}\in\mathbb{R}^{m^{2}}$
denote the respective load vectors. Moreover, the interface condition
can also be formulated as the linear equation $\boldsymbol{C}_{1}^{\mathbf{T}}\boldsymbol{u}_{1}^{\prime}=\boldsymbol{C}_{2}^{\mathbf{T}}\boldsymbol{u}_{2}^{\prime}$,
where $\boldsymbol{C}_{1},\boldsymbol{C}_{2}\in\mathbb{R}^{m^{2}\times m}$
are the interface matrices. 

Following the finite element tearing
and interconnect (FETI) approach [\hyperlink{ref35}{35}], let $\boldsymbol{\lambda}\in\mathbb{R}^{m}$
denote the Lagrange multiplier of the interface equation, $\boldsymbol{u}_{1}=\left[\boldsymbol{u}_{1}^{\prime};\boldsymbol{\lambda}\right]\in\mathbb{R}^{n_{1}}\equiv\mathbb{R}^{m^{2}+m}$,
$\boldsymbol{u}_{2}=\left[\boldsymbol{u}_{2}^{\prime}\right]\in\mathbb{R}^{n_{2}}\equiv\mathbb{R}^{m^{2}}$
denote the solution variables, and $\boldsymbol{v}_{1}=\boldsymbol{\lambda}\left(\boldsymbol{u}_{1}\right)\in\mathbb{R}^{m_{2}}\equiv\mathbb{R}^{m}$,
$\boldsymbol{v}_{2}=\boldsymbol{C}_{2}^{\mathbf{T}}\boldsymbol{u}_{2}\in\mathbb{R}^{m_{2}}\equiv\mathbb{R}^{m}$
denote the coupling variables. Subsequently, as per \hyperlink{eq27}{Eq. 2.7}, we can
formulate a modular algebraic system with following component residuals
and interface functions. 
\begin{align}
\boldsymbol{f}_{1}\left(\boldsymbol{u}_{1};\boldsymbol{v}_{2},\boldsymbol{\xi}_{1}\right) & =\left[\begin{array}{cc}
\boldsymbol{A}_{1}\left(\boldsymbol{\xi}_{1}\right) & \boldsymbol{C}_{1}\\
\boldsymbol{C}_{1}^{\mathbf{T}} & 0
\end{array}\right]\boldsymbol{u}_{1}-\left[\begin{array}{c}
\boldsymbol{b}_{1}\\
\boldsymbol{v}_{2}
\end{array}\right], & \boldsymbol{g}_{1}\left(\boldsymbol{u}_{1}\right)=\boldsymbol{\lambda}\left(\boldsymbol{u}_{1}\right),\nonumber \\
\boldsymbol{f}_{2}\left(\boldsymbol{u}_{2};\boldsymbol{v}_{1},\boldsymbol{\xi}_{2}\right) & =\boldsymbol{A}_{2}\left(\boldsymbol{\xi}_{2}\right)\boldsymbol{u}_{2}-\boldsymbol{C}_{2}\boldsymbol{v}_{1}-\boldsymbol{b}_{2}, & \boldsymbol{g}_{2}\left(\boldsymbol{u}_{2}\right)=\boldsymbol{C}_{2}^{\mathbf{T}}\boldsymbol{u}_{2}.
\end{align}

The probability density function of the total energy $E:\forall\boldsymbol{\xi}\in\Xi$, 
\begin{equation}
E\left(\boldsymbol{\xi}\right)=\frac{1}{2}\left(\int_{\Omega_{1}}u_{1}\left(\boldsymbol{x},\boldsymbol{\xi}\right)^{2}d\boldsymbol{x}+\int_{\Omega_{2}}u_{2}\left(\boldsymbol{x},\boldsymbol{\xi}\right)^{2}d\boldsymbol{x}\right)
\end{equation}
and the statistics of $u_{1},u_{2}$ are the quantities of interest
in this study. The numerical values of the deterministic parameters
used in this study are listed in \hyperlink{tab1}{Table 1}.

\begin{table}[htbp]
\hypertarget{tab1}{}
\caption{Deterministic parameter values in the Poisson problem.}
\begin{center}\scriptsize
\renewcommand{\arraystretch}{1.3}
\begin{tabular}{cccccccc}
\toprule 
$\bar{a}_{1}$ & $\bar{a}_{2}$ & $b_{1}$ & $b_{2}$ & $\delta_{1}$ & $\delta_{2}$ & $l_{1}$ & $l_{2}$\tabularnewline
\midrule
\midrule 
$0.5$ & $1.0$ & $4.0$ & $-4.0$ & $0.5$ & $0.2$ & $0.2$ & $0.5$\tabularnewline
\bottomrule
\end{tabular}
\par\end{center}
\end{table}

\subsubsection{Modular deterministic solver: Setup and verification}

The deterministic solver component $\boldsymbol{m}_{1}$ and $\boldsymbol{m}_{2}$ wer developed as $\mathtt{MATLAB}\text{\texttrademark}$
function modules, in which the output respective solution updates are
computed using Newton's method as follows. 
\begin{align}
\boldsymbol{m}_{1}\left(\boldsymbol{v}_{2},\boldsymbol{\xi}_{1}\right)  =\left(\frac{\partial\boldsymbol{f}_{1}}{\partial\boldsymbol{u}_{1}}\left(\boldsymbol{\xi}_{1}\right)\right)^{-1}\left[\begin{array}{c}
\boldsymbol{b}_{1}\\
\boldsymbol{v}_{2}
\end{array}\right],\nonumber \\
\boldsymbol{m}_{2}\left(\boldsymbol{v}_{1},\boldsymbol{\xi}_{2}\right) =\left(\frac{\partial\boldsymbol{f}_{2}}{\partial\boldsymbol{u}_{2}}\left(\boldsymbol{\xi}_{2}\right)\right)^{-1}\left(\boldsymbol{C}_{2}\boldsymbol{v}_{1}+\boldsymbol{b}_{2}\right).
\end{align}

Moreover, to accelerate convergence, a modified relaxed BGS approach
was implemented with the optimal value of $0.9$ for both relaxation
factors. Subsequently, the second order accuracy of the
solver was verified using the method of manufactured solutions (MMS) [\hyperlink{ref36}{36}].
Further details are provided in \hyperlink{appB1}{Appendix B}.

\subsubsection{NISP based uncertainty propagation}

Following the verification study, both the standard and reduced NISP
based uncertainty propagation algorithms were implemented by
reusing deterministic solver components $\boldsymbol{m}_{1}$ and $\boldsymbol{m}_{2}$. Tolerance values
of $\epsilon_{1,\mathrm{dim}}=10^{-4}$, $\epsilon_{2,\mathrm{dim}}=10^{-5}$,
$\epsilon_{1,\mathrm{ord}}=\epsilon_{2,\mathrm{ord}}=10^{-3}$
were used in the reduced NISP method implementation. For $m=31$, $s_{1}=s_{2}=4$,
$p=4$ , the probability density function of $E$ and the first two solution moments were
computed using the converged gPC coefficient matrices $\hat{\boldsymbol{U}}_{1}$
and $\hat{\boldsymbol{U}}_{2}$ from both algorithms. The results are
shown \hyperlink{fig2}{Figure 2} and \hyperlink{fig3}{Figure 3} respectively. 

\begin{figure}
    \hypertarget{fig2}{}
    \centering
    \includegraphics[bb=40bp 220bp 552bp 589bp,clip,scale=0.55]{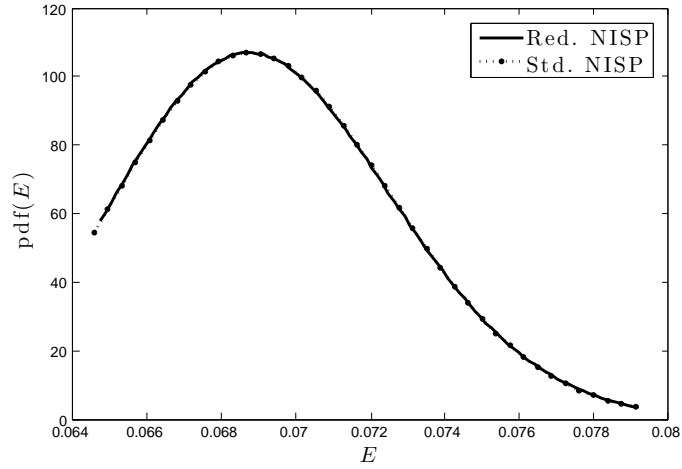} 
    \caption{Probability density function of the total energy $E$ computed
using both NISP method implementations. The densities were computed using the KDE method with $10^{5}$ samples.}
    \label{fig:fig2}
\end{figure}

\begin{figure}
    \hypertarget{fig3}{}
    \centering
    \includegraphics[bb=100bp 350bp 1100bp 800bp,scale=0.50]{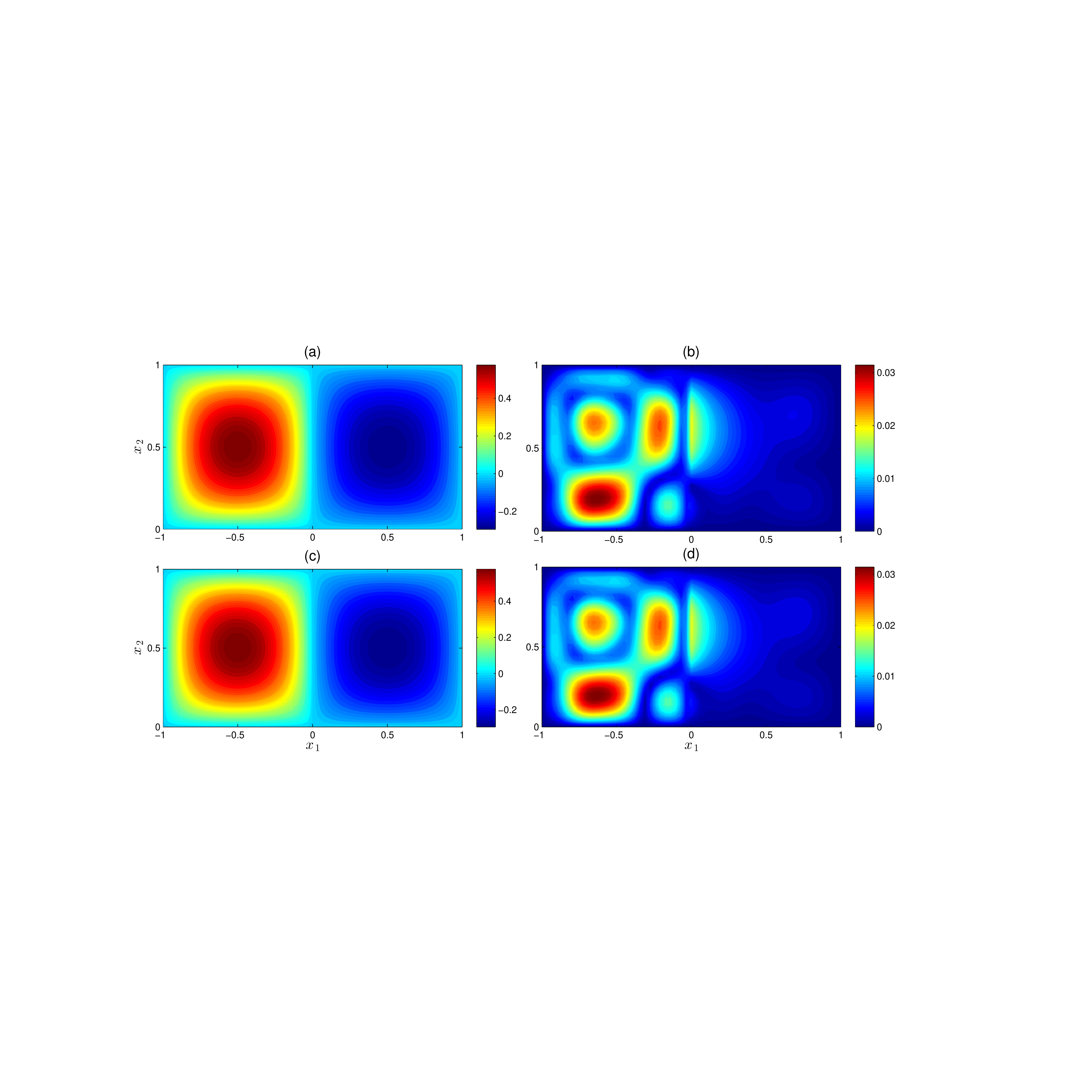}
    \caption{First two moments of the solution computed
using both NISP method implementations. Subfigures (a, c) and (b, d) correspond to the mean and standard deviation respectively, while (a, b) and (c, d) correspond to the reduced and standard NISP method implementations respectively.}
    \label{fig:fig3}
\end{figure}

For both algorithms, keeping the mesh size and tolerances values unchanged, the
approximation errors $\varepsilon_{s},\varepsilon_{r}$ and computational
costs (wall-times) $\mathcal{C}_{s},\mathcal{C}_{r}$ for several instances of $s_{1}$,
$s_{2}$ and $p$ are listed in \hyperlink{tab2}{Table 2}. Moreover, the parameters
of the reduced NISP algorithm in last iteration are also listed. The
number of iterations in all tests were observed to be between $9$
and $11$ in each implementation, which implies that the convergence
rate in both algorithms are more or less invariant with respect to the stochastic dimensions $s_{1}$, $s_{2}$ and the order of accuracy $p$. While the choice of the tolerances in the reduced NISP algorithm seem
arbitrary, they were in fact selected based on the results of preliminary
tests using a coarser mesh. \hyperlink{fig4}{Figure 4} illustratess the effects of varying $\epsilon_{1,\mathrm{dim}}$
and $\epsilon_{2,\mathrm{dim}}$ on the standard deviation of the solution, for $m=11$ and all other parameters unchanged.

The highest speedup factor observed 
 is $\approx 11.1$.  Moreover, the approximation errors in the standard NISP method implementation
are observed to decay exponentially, which indicates a high degree
of regularity in the stochastic solutions. Furthermore, for the reduced NISP method implementation, the asymptotic upper bound 
predicted in \hyperlink{sec34}{\S 3.4}, is observed in its approximation
errors.

\begin{table}[htbp]
\hypertarget{tab2}{}
\caption{Comparison of the approximation errors and computational costs
(seconds) obtained in the standard and reduced NISP method implementations for the Poisson problem. The
reduced dimensions and orders are also listed here.}
\begin{center}\scriptsize
\renewcommand{\arraystretch}{1.3}
\begin{tabular}{cc|ccc|cccccccc|c}
\toprule 
 &  & \multicolumn{3}{c}{Standard NISP} \vline & \multicolumn{8}{c}{Reduced NISP}\vline & \tabularnewline
{$s_{1}$, $s_{2}$} & {$p$} & {$Q$} & {$\varepsilon_{s}$} & {$\mathcal{C}_{s}$ } & {$d_{1}$} & {$\tilde{p}_{1}$} & {$\tilde{Q}_{1}$} & {$d_{2}$} & {$\tilde{p}_{2}$} & {$\tilde{Q}_{2}$} & {$\varepsilon_{r}$} & {$\mathcal{C}_{r}$} & {$\mathcal{C}_{s}/ \mathcal{C}_{r}$} \tabularnewline
\midrule
\midrule 
 & {$2$} & {$85$} & {$2.9\times10^{-3}$} & {$6$} & {$3$} & {$1$} & {$35$} & {$5$} & {$1$} & {$85$} & {$1.3\times10^{-2}$} & {$5$} & {$1.2$}\tabularnewline
\cmidrule{2-14} 
{$3$} & {$3$} & {$389$} & {$9.2\times10^{-4}$} & {$32$} & {$3$} & {$2$} & {$84$} & {$5$} & {$1$} & {$126$} & {$6.6\times10^{-3}$} & {$8$} & {$4.0$}\tabularnewline
\cmidrule{2-14}  
 & {$4$} & {$1457$} & {$4.4\times10^{-4}$} & {$143$} & {$3$} & {$3$} & {$165$} & {$5$} & {$1$} & {$126$} & {$4.4\times10^{-3}$} & {$13$} & {$11.0$}\tabularnewline
\midrule 
 & {$2$} & {$145$} & {$3.2\times10^{-3}$} & {$16$} & {$4$} & {$1$} & {$70$} & {$6$} & {$1$} & {$145$} & {$1.3\times10^{-2}$} & {$11$} & {$1.5$}\tabularnewline
\cmidrule{2-14}  
{$4$} & {$3$} & {$849$} & {$1.0\times10^{-3}$} & {$79$} & {$4$} & {$2$} & {$210$} & {$6$} & {$1$} & {$210$} & {$6.5\times10^{-3}$} & {$21$} & {$3.8$}\tabularnewline
\cmidrule{2-14}  
 & {$4$} & {$3937$} & {$4.6\times10^{-4}$} & {$544$} & {$4$} & {$3$} & {$495$} & {$6$} & {$1$} & {$210$} & {$4.4\times10^{-3}$} & {$45$} & {$12.1$}\tabularnewline
\midrule 
 & {$2$} & {$221$} & {$3.1\times10^{-3}$} & {$17$} & {$5$} & {$1$} & {$126$} & {$7$} & {$1$} & {$221$} & {$1.5\times10^{-2}$} & {$25$} & {$0.7$}\tabularnewline
\cmidrule{2-14}  
{$5$} & {$3$} & {$1581$} & {$1.1\times10^{-3}$} & {$176$} & {$5$} & {$2$} & {$462$} & {$7$} & {$1$} & {$330$} & {$8.7\times10^{-3}$} & {$47$} & {$3.7$}\tabularnewline
\cmidrule{2-14}  
 & {$4$} & {$8801$} & {$4.7\times10^{-4}$} & {$1751$} & {$5$} & {$3$} & {$1287$} & {$7$} & {$2$} & {$1716$} & {$8.3\times10^{-3}$} & {$509$} & {$3.4$}\tabularnewline
\bottomrule
\end{tabular}
\par\end{center}
\end{table}.

\begin{figure}
    \hypertarget{fig4}{}
    \centering
    \includegraphics[bb=90bp 350bp 1100bp 800bp,scale=0.50]{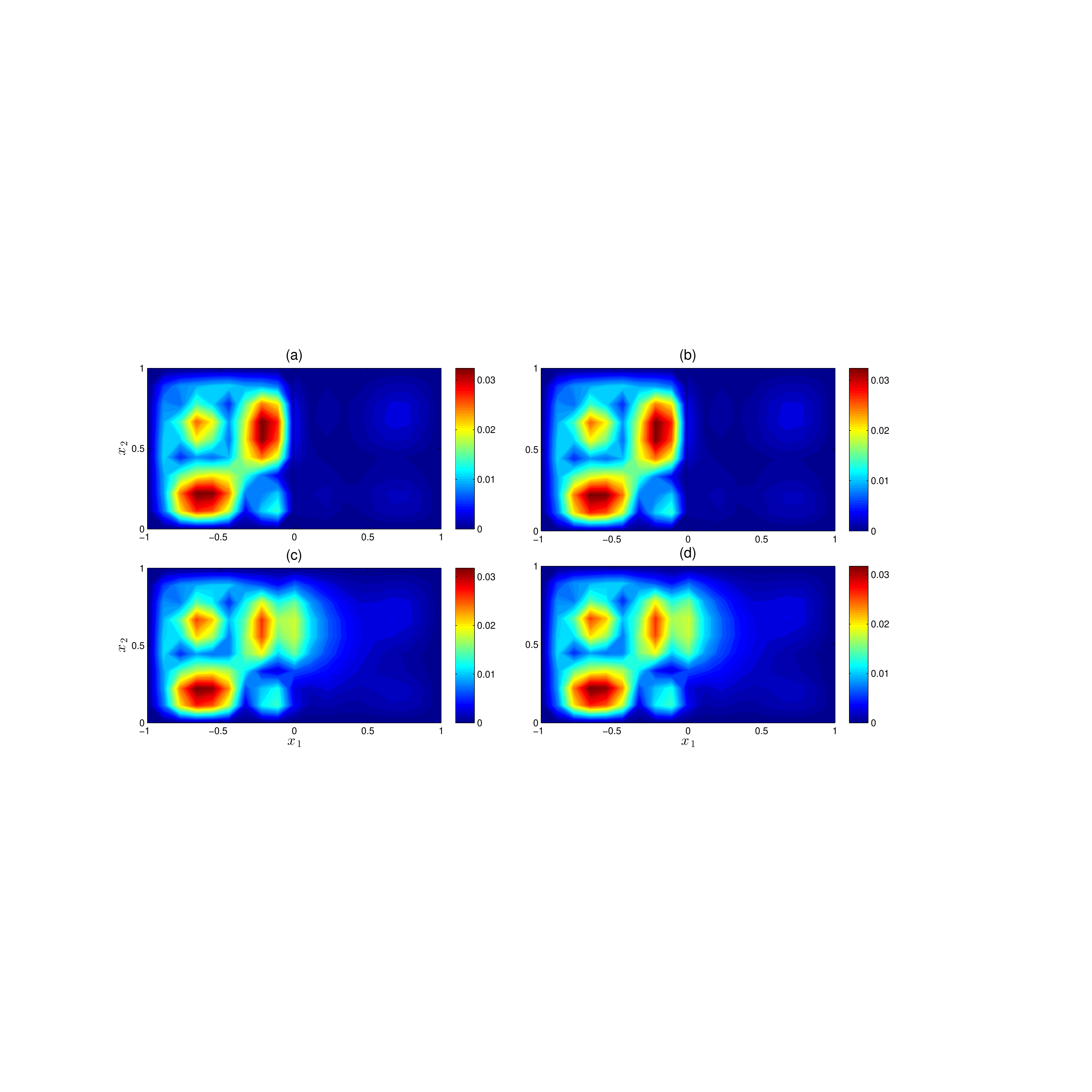}
    \caption{Comparison of computed standard deviations on coarser mesh. Subfigures (a, b, c) correspond to the reduced NISP method implementation while (d) corresponds to the standard NISP method implementation. The tolerance values are as follows. (a) $\epsilon_{1,\mathrm{dim}}=10^{-4}$, $\epsilon_{2,\mathrm{dim}}=10^{-4}$, (b) $\epsilon_{1,\mathrm{dim}}=10^{-5}$, $\epsilon_{2,\mathrm{dim}}=10^{-4}$, (c) $\epsilon_{1,\mathrm{dim}}=10^{-4}$, $\epsilon_{2,\mathrm{dim}}=10^{-5}$.}
    \label{fig:fig4}
\end{figure}

\subsection{Boussinesq flow problem}
\hypertarget{sec42}{}
\begin{figure}
    \hypertarget{fig5}{}
    \centering
    \includegraphics[bb=125bp 500bp 475bp 750bp,scale=0.65]{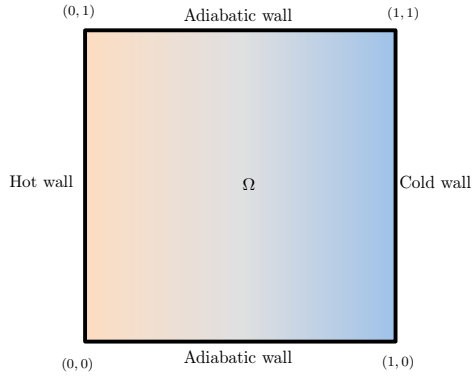}
    \caption{Computational domain for the Boussinesq flow problem.}
    \label{fig:fig5}
\end{figure}

The Boussinesq model describes thermally driven, incompressible
flows and is widely used in oceanic and atmospheric modeling [\hyperlink{ref37}{37}, \hyperlink{ref38}{38}]. Here, we consider a multi-physics setup with uncertain fluid properties and boundary conditions.

\subsubsection{Model setup}

Let $\Omega\equiv\left(0,1\right)_{x_{1}}\times\left(0,1\right)_{x_{2}}$
denote the spatial domain and \textbf{$\boldsymbol{u}=\left[\begin{array}{cc}
u_{1} & u_{2}\end{array}\right]^{\mathbf{T}}$} , $p$ $,T$ denote the non-dimensional [\hyperlink{ref39}{39}] fluid velocity,
pressure and temperature respectively. The governing equations for
the fluid variables are as follows. $\forall \boldsymbol{\xi} \in \Xi$,

\begin{align}
\boldsymbol{\nabla}^{\mathbf{T}}\boldsymbol{u}\left(\boldsymbol{x},\boldsymbol{\xi}\right) & =0,\nonumber \\
\left(\boldsymbol{u}\left(\boldsymbol{x},\boldsymbol{\xi}\right)^{\mathbf{T}}\boldsymbol{\nabla}\right)\boldsymbol{u}\left(\boldsymbol{x},\boldsymbol{\xi}\right)+\boldsymbol{\nabla}p\left(\boldsymbol{x},\boldsymbol{\xi}\right)\nonumber \\
-\mathrm{Pr}\boldsymbol{\nabla}^{\mathbf{T}}\boldsymbol{\nabla}\boldsymbol{u}\left(\boldsymbol{x},\boldsymbol{\xi}\right)-\mathrm{Pr}\mathrm{Ra}\left(\boldsymbol{x},\boldsymbol{\xi}_{1}\left(\boldsymbol{\xi}\right)\right)T\left(\boldsymbol{x},\boldsymbol{\xi}\right)\boldsymbol{e}_{2} & =\boldsymbol{0},\nonumber \\
\left(\boldsymbol{u}\left(\boldsymbol{x},\boldsymbol{\xi}\right)^{\mathbf{T}}\boldsymbol{\nabla}\right)T\left(\boldsymbol{x},\boldsymbol{\xi}\right)-\boldsymbol{\nabla}^{\mathbf{T}}\boldsymbol{\nabla}T\left(\boldsymbol{x},\boldsymbol{\xi}\right) & =0, & \boldsymbol{x}\in\Omega,
\end{align}
with homogenous Dirichlet boundary conditions for $\boldsymbol{u}$
and Neumann boundary conditions for $p$ at all boundaries. The boundary conditions for temperature, as shown in \hyperlink{fig5}{Figure 5}, are as follows.
\begin{align}
\frac{\partial T}{\partial x_{2}}\left(x_{1},0,\boldsymbol{\xi}\right)=\frac{\partial T}{\partial x_{2}}\left(x_{1},1,\boldsymbol{\xi}\right) & =0, & x_{1}\in\left[0,1\right],\nonumber \\
T\left(0,x_{2},\boldsymbol{\xi}\right)-T_{h}\left(x_{2},\boldsymbol{\xi}_{2}\right)=T\left(1,x_{2},\boldsymbol{\xi}\right) & =0, & x_{2}\in\left[0,1\right].
\end{align}

Here, $\boldsymbol{e}_{2}$ denotes $\left[\begin{array}{cc}
0 & 1\end{array}\right]^{\mathbf{T}}$ while $\mathrm{Pr}$ and $\mathrm{Ra}$ denote the Prandtl and Rayleigh
numbers respectively. $T_{h}$ denotes the hot-wall temperature such
that $\forall x_{2}\in\left[0,1\right],\boldsymbol{\xi}_{2}\in\Xi_2$,
\begin{equation}
T_{h}\left(x_{2},\boldsymbol{\xi}_{2}\right)=\bar{T}_{h}+h\left(x_{2},\boldsymbol{\xi}_{2}\right)\sin^{2}\left(\pi x_{2}\right),
\end{equation}
where $\bar{T}_{h}$ is the mean hot-wall temperature and $h$ denotes
the perturbation amplitude. 

In this study, $\mathrm{Ra}$ and $h$ are assumed to be independent
random fields, and modeled using the following KL expansions. $\forall\boldsymbol{x}\in\Omega,\boldsymbol{\xi}_{1}\in\Xi_{1}$,
\begin{equation}
\mathrm{Ra}\left(\boldsymbol{x},\boldsymbol{\xi}_{1}\right)=\bar{\mathrm{Ra}}+\sqrt{3}\delta_{\mathrm{Ra}}\sum_{j=1}^{s_{1}}\gamma_{\mathrm{Ra},j}\left(\boldsymbol{x}\right)\xi_{1j},
\end{equation}
where $\bar{\mathrm{Ra}}$ denotes the mean of $\mathrm{Ra}$ and
$\left\{ \xi_{1j}\sim U\left[-1,1\right]\right\}_{j=1}^{s_{1}} $
are $i.i.d.$ random variables. Similarly, $\forall x_{2}\in\left(0,1\right),\boldsymbol{\xi}_{2}\in\Xi_{2}$,
\begin{equation}
h\left(x_{2},\boldsymbol{\xi}_{2}\right)=\sqrt{3}\delta_{h}\sum_{j=1}^{s_{2}}\gamma_{h,j}\left(\boldsymbol{x}\right)\xi_{2j},
\end{equation}
where $\left\{ \xi_{2j}\sim U\left[-1,1\right]\right\}_{j=1}^{s_{2}} $
are $i.i.d.$ random variables. Moreover, we assume that both $\mathrm{Ra}$
and $h$ have exponential covariance kernels
\begin{align}
C_{\mathrm{Ra}}\left(\boldsymbol{x},\boldsymbol{y}\right) & =\delta_{\mathrm{Ra}}^{2}\exp\left(-\frac{\left\Vert \boldsymbol{x}-\boldsymbol{y}\right\Vert _{1}}{l_{\mathrm{Ra}}}\right),&\boldsymbol{x},\boldsymbol{y}\in\Omega,\nonumber \\
C_{h}\left(x_{2},y_{2}\right) & =\delta_{h}^{2}\exp\left(-\frac{\left|x_{2}-y_{2}\right|}{l_{h}}\right),&x_{2},y_{2}\in\left[0,1\right],
\end{align}
where $\delta_{\mathrm{Ra}}$, $\delta_{h}$ denote the respective
coefficients of variations, and $l_{\mathrm{Ra}}$, $l_{h}$ denote
the respective correlation lengths. The analytic expressions for $\left\{\gamma_{\mathrm{Ra},j}\right\}_{j>0}$ and $\left\{\gamma_{h,j}\right\}_{j>0}$
are provided in \hyperlink{appA}{Appendix A}. Furthermore, the pressure Poisson equation
\begin{equation}
\boldsymbol{\nabla}^{\mathbf{T}}\boldsymbol{\nabla}p\left(\boldsymbol{x},\boldsymbol{\xi}\right)+\boldsymbol{\nabla}^{\mathbf{T}}\left(\left(\boldsymbol{u}\left(\boldsymbol{x},\boldsymbol{\xi}\right)^{\mathbf{T}}\boldsymbol{\nabla}\right)\boldsymbol{u}\left(\boldsymbol{x},\boldsymbol{\xi}\right)-\mathrm{Pr}\mathrm{Ra}\left(\boldsymbol{x},\boldsymbol{\xi}_{1}\right)T\left(\boldsymbol{x},\boldsymbol{\xi}\right)\boldsymbol{e}_{2}\right)=0
\end{equation}
is used in place of the continuity equation, to close the momentum component of the PDE system. 

Each component PDE system is spatially discretized using a
finite volume method, with linear central-differencing schemes [\hyperlink{ref40}{40}], on a uniform
grid with $m\times m$ cells. Let $\boldsymbol{u}_{1}^{\prime},\boldsymbol{u}_{2}^{\prime},\boldsymbol{p}^{\prime},\boldsymbol{t}^{\prime}\in\mathbb{R}^{m^{2}}$
denote the respective vectors of cell-centroidal horizontal velocity,
vertical velocity, pressure and temperature, which solve the nonlinear
system
\hypertarget{eq415}{}
\begin{align}
\left(\boldsymbol{K}_{u}+\boldsymbol{A}\left(\boldsymbol{u}_{1}^{\prime},\boldsymbol{u}_{2}^{\prime}\right)\right)\boldsymbol{u}_{1}^{\prime}+\boldsymbol{B}_{1}\boldsymbol{p}^{\prime} & =\boldsymbol{0},\nonumber \\
\left(\boldsymbol{K}_{u}+\boldsymbol{A}\left(\boldsymbol{u}_{1}^{\prime},\boldsymbol{u}_{2}^{\prime}\right)\right)\boldsymbol{u}_{2}^{\prime}+\boldsymbol{B}_{2}\boldsymbol{p}^{\prime}-\boldsymbol{R}\left(\boldsymbol{\xi}_{1}\right)\boldsymbol{t}^{\prime} & =\boldsymbol{0},\nonumber \\
\boldsymbol{K}_{p}\boldsymbol{p}^{\prime}+\boldsymbol{C}_{1}\left(\boldsymbol{u}_{1}^{\prime},\boldsymbol{u}_{2}^{\prime}\right)\boldsymbol{u}_{1}^{\prime}+\boldsymbol{C}_{2}\left(\boldsymbol{u}_{1}^{\prime},\boldsymbol{u}_{2}^{\prime}\right)\boldsymbol{u}_{2}^{\prime}-\boldsymbol{S}\left(\boldsymbol{\xi}_{1}\right)\boldsymbol{t}^{\prime} & =\boldsymbol{0},\nonumber \\
\left(\boldsymbol{K}_{T}+\boldsymbol{A}\left(\boldsymbol{u}_{1}^{\prime},\boldsymbol{u}_{2}^{\prime}\right)\right)\boldsymbol{t}^{\prime}-\boldsymbol{h}\left(\boldsymbol{\xi}_{2}\right) & =\boldsymbol{0}.
\end{align}
where each term in \hyperlink{eq415}{Eq. 4.15} denotes its respective discretized operator
in the coupled PDE system. Subsequently, we formulate a modular multi-physics
setup by separating the momentum and energy components of the coupled
algebraic system. As per \hyperlink{eq25}{Eq. 2.5}, let $\boldsymbol{u}_{1}=\left[\boldsymbol{u}_{1}^{\prime};\boldsymbol{u}_{2}^{\prime};\boldsymbol{p}^{\prime}\right]\in\mathbb{R}^{n_{1}}\equiv\mathbb{R}^{3m^{2}}$,
$\boldsymbol{u}_{2}=\boldsymbol{t}^{\prime}\in\mathbb{R}^{n_{2}}\equiv\mathbb{R}^{m^{2}}$
denote the respective solution variables in the modular algebraic
system. The component residuals are defined as follows.
\begin{align}
\boldsymbol{f}_{1}\left(\boldsymbol{u}_{1};\boldsymbol{u}_{2},\boldsymbol{\xi}_{1}\right) & =\left[\begin{array}{cc}
\boldsymbol{K}_{u}+\boldsymbol{A}\left(\boldsymbol{u}_{1}^{\prime}\left(\boldsymbol{u}_{1}\right),\boldsymbol{u}_{2}^{\prime}\left(\boldsymbol{u}_{1}\right)\right) & \boldsymbol{0}\\
\boldsymbol{0} & \boldsymbol{K}_{u}+\boldsymbol{A}\left(\boldsymbol{u}_{1}^{\prime}\left(\boldsymbol{u}_{1}\right),\boldsymbol{u}_{2}^{\prime}\left(\boldsymbol{u}_{1}\right)\right)\\
\boldsymbol{C}_{1}\left(\boldsymbol{u}_{1}^{\prime}\left(\boldsymbol{u}_{1}\right),\boldsymbol{u}_{2}^{\prime}\left(\boldsymbol{u}_{1}\right)\right) & \boldsymbol{C}_{2}\left(\boldsymbol{u}_{1}^{\prime}\left(\boldsymbol{u}_{1}\right),\boldsymbol{u}_{2}^{\prime}\left(\boldsymbol{u}_{1}\right)\right)
\end{array}\right.\nonumber \\
 & \left.\begin{array}{c}
\boldsymbol{B}_{1}\\
\boldsymbol{B}_{2}\\
\boldsymbol{K}_{p}
\end{array}\right]\boldsymbol{u}_{1}-\left[\begin{array}{c}
\boldsymbol{0}\\
\boldsymbol{R}\left(\boldsymbol{\xi}_{1}\right)\\
\boldsymbol{S}\left(\boldsymbol{\xi}_{1}\right)
\end{array}\right]\boldsymbol{u}_{2},\nonumber \\
\boldsymbol{f}_{2}\left(\boldsymbol{u}_{2};\boldsymbol{u}_{1},\boldsymbol{\xi}_{2}\right) & =\left(\boldsymbol{K}_{T}+\boldsymbol{A}\left(\boldsymbol{u}_{1}^{\prime}\left(\boldsymbol{u}_{1}\right),\boldsymbol{u}_{2}^{\prime}\left(\boldsymbol{u}_{1}\right)\right)\right)\boldsymbol{u}_{2}-\boldsymbol{h}\left(\boldsymbol{\xi}_{2}\right).
\end{align}

The quantities of interest in this study are the probability density functions of the
(scaled) kinetic energy $K$ and thermal energy $E$: $\forall\boldsymbol{\xi}\in\Xi$,
\begin{equation}
K\left(\boldsymbol{\xi}\right)=\frac{1}{2}\left(\int_{\Omega}u_{1}\left(\boldsymbol{x},\boldsymbol{\xi}\right)^{2}d\boldsymbol{x}+\int_{\Omega}u_{2}\left(\boldsymbol{x},\boldsymbol{\xi}\right)^{2}d\boldsymbol{x}\right),\ E\left(\boldsymbol{\xi}\right)=\int_{\Omega}T\left(\boldsymbol{x},\boldsymbol{\xi}\right)d\boldsymbol{x},
\end{equation}
and the statistics of the fluid velocity and temperature. \hyperlink{tab3}{Table 3}
lists the numerical values of the deterministic parameters used in
this study.

\begin{table}[htbp]
\hypertarget{tab3}{}
\caption{Deterministic parameter values in the Boussinesq flow problem.}
\begin{center}\scriptsize
\renewcommand{\arraystretch}{1.3}
\begin{tabular}{ccccccc}
\toprule 
$\mathrm{Pr}$ & $\bar{\mathrm{Ra}}$ & $\bar{T}_{h}$ & $\delta_{\mathrm{Ra}}$ & $\delta_{h}$ & $l_{\mathrm{Ra}}$ & $l_{h}$\tabularnewline
\midrule
\midrule 
$0.71$ & $1000$ & $1$ & $200$ & $0.5$ & $0.5$ & $0.5$\tabularnewline
\bottomrule
\end{tabular}
\par\end{center}
\end{table}

\subsubsection{Modular deterministic solver: Setup and verification}

The solver components $\boldsymbol{m}_{1}$ and $\boldsymbol{m}_{2}$ were developed as
$\mathtt{MATLAB}\text{\texttrademark}$ function modules, with each module solving
a Newton system to compute its respective solution updates. Therefore,
\begin{align}
\boldsymbol{m}_{1}\left(\boldsymbol{u}_{1},\boldsymbol{u}_{2},\boldsymbol{\xi}_{1}\right) & =\boldsymbol{u}_{1}-\left(\frac{\partial\boldsymbol{f}_{1}}{\partial\boldsymbol{u}_{1}}\left(\boldsymbol{u}_{1};\boldsymbol{u}_{2},\boldsymbol{\xi}_{1}\right)\right)^{-1}\boldsymbol{f}_{1}\left(\boldsymbol{u}_{1};\boldsymbol{u}_{2},\boldsymbol{\xi}_{1}\right),\nonumber \\
\boldsymbol{m}_{2}\left(\boldsymbol{u}_{2},\boldsymbol{u}_{1},\boldsymbol{\xi}_{2}\right) & =\boldsymbol{u}_{2}-\left(\frac{\partial\boldsymbol{f}_{2}}{\partial\boldsymbol{u}_{2}}\left(\boldsymbol{u}_{2};\boldsymbol{u}_{1},\boldsymbol{\xi}_{2}\right)\right)^{-1}\boldsymbol{f}_{2}\left(\boldsymbol{u}_{2};\boldsymbol{u}_{1},\boldsymbol{\xi}_{2}\right).
\end{align}

Using MMS, a verification study was carried out and second order accuracy was observed in the numerical solver. The details are provided in \hyperlink{appB2}{Appendix B}.

\subsubsection{NISP based uncertainty propagation}

We implemented both NISP based uncertainty
propagation methods by reusing the deterministic solver components $\boldsymbol{m}_{1}$ and $\boldsymbol{m}_{2}$. The tolerance
values used in the reduced NISP method implementation are $\epsilon_{1,\mathrm{dim}}=10^{-3},\epsilon_{2,\mathrm{dim}}=10^{-4}$
and $\epsilon_{1,\mathrm{ord}}=\epsilon_{2,\mathrm{ord}}=10^{-3}$.
Subsequently, for $m=25$, $s_{1}=s_{2}=3$ and $p=4$, the probability density function of
$K$ and $E$, and the first two moments of $u_{1},u_{2},T$, were computed and compared. The results are shown in 
\hyperlink{fig6}{Figure 6}, \hyperlink{fig7}{Figure 7} and \hyperlink{fig8}{Figure 8} respectively.  

\begin{figure}
    \hypertarget{fig6}{}
    \centering
    \includegraphics[bb=160bp 400bp 1050bp 750bp,clip,scale=0.55]{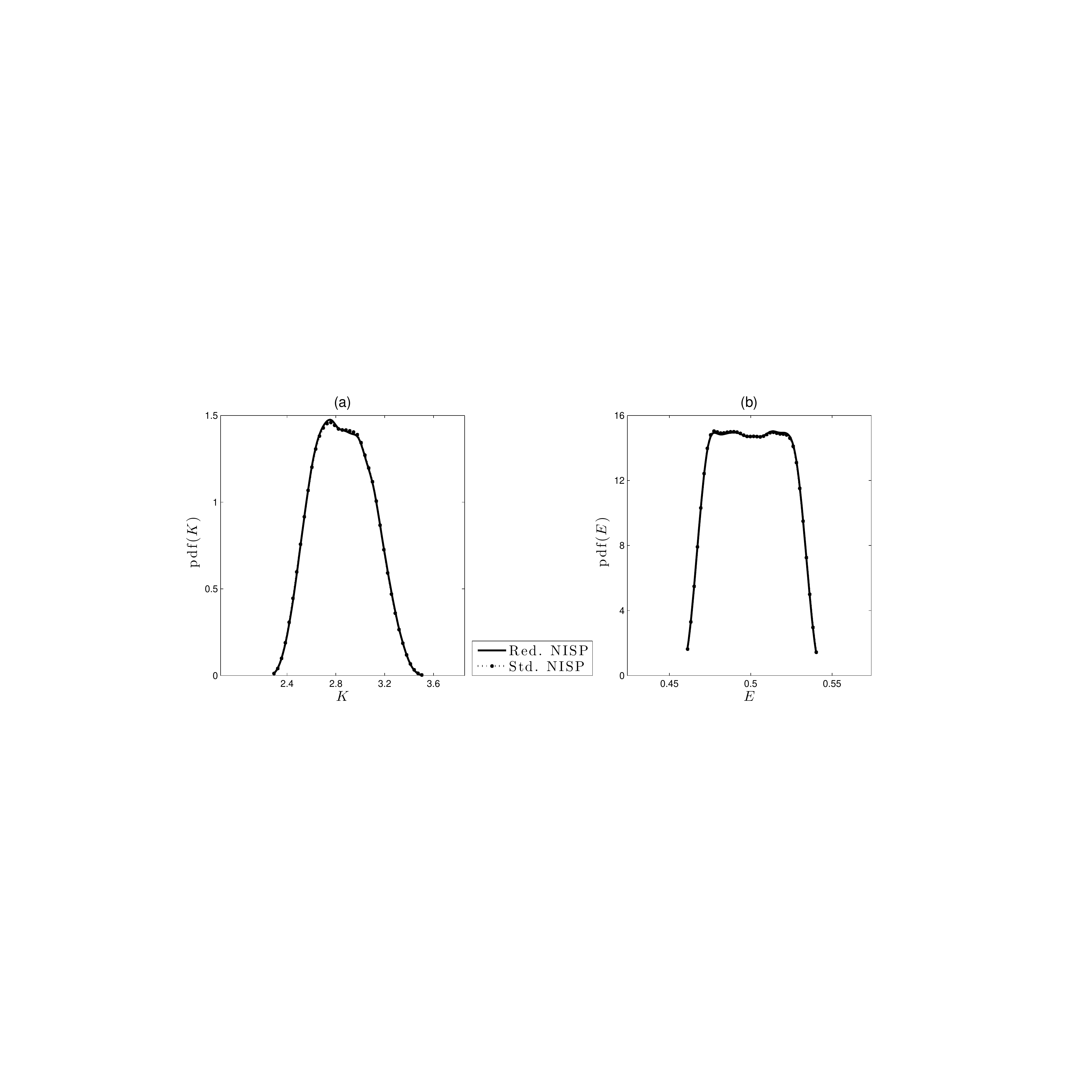} 
    \caption{Probability density function of the fluid energies computed
using both NISP method implementations. Subfigure (a) corresponds to the kinetic energy $K$, while (b) corresponds to the thermal energy $E$. The densities were computed using the KDE method with $10^{5}$ samples.}
    \label{fig:fig6}
\end{figure}

\begin{figure}
    \hypertarget{fig7}{}
    \centering
    \includegraphics[bb=80bp 340bp 1100bp 820bp,scale=0.5]{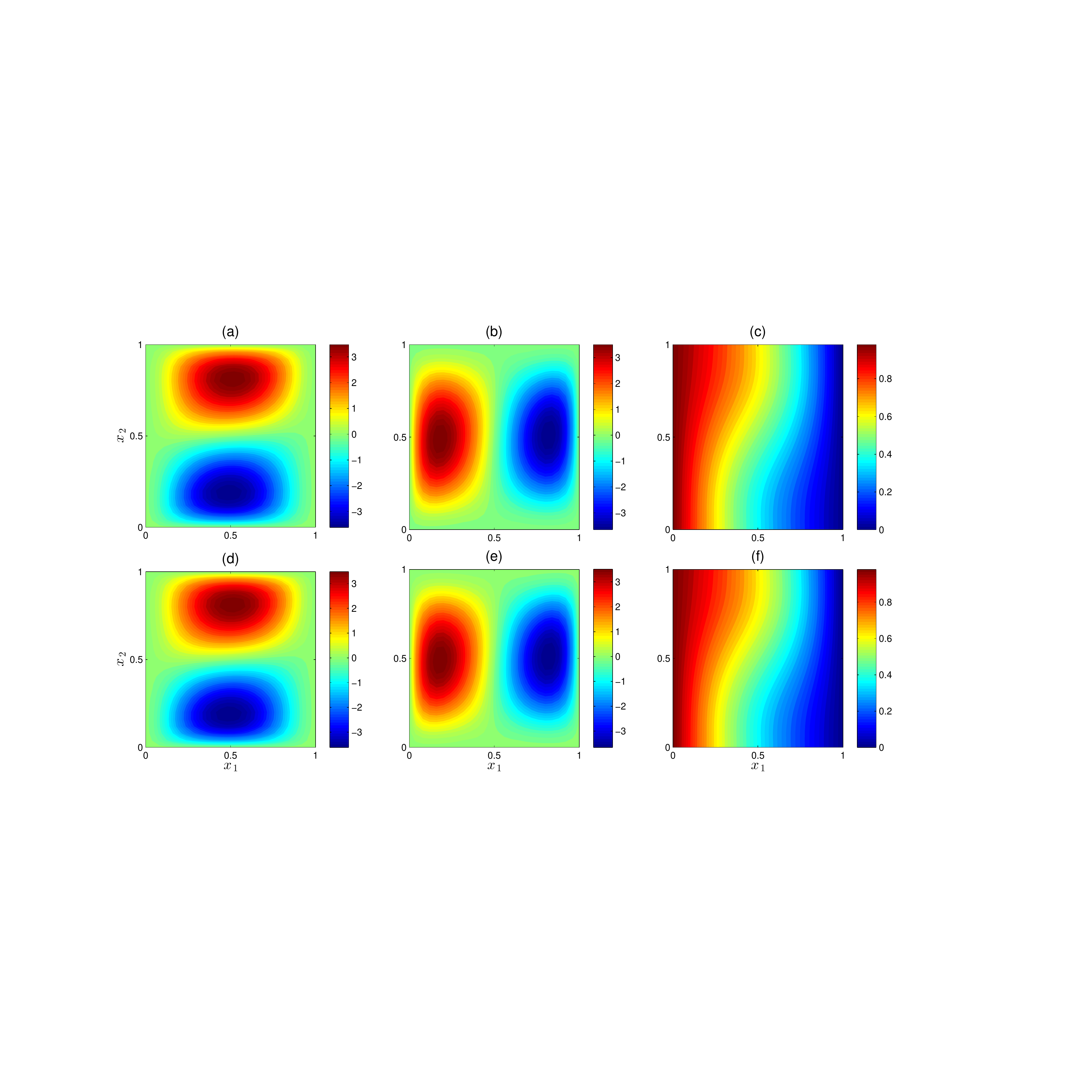}
    \caption{Mean of solution fields obtained using both NISP method implementations. Subfigures (a, d), (b, e) and (c, f) correspond to $u_1$, $u_2$ and $T$ respectively, while (a, b, c) and (d, e, f) correspond to the reduced and standard NISP method implementations respectively.}
    \label{fig:fig7}
\end{figure}

\begin{figure}
    \hypertarget{fig8}{}
    \centering
    \includegraphics[bb=115bp 340bp 1100bp 800bp,scale=0.54]{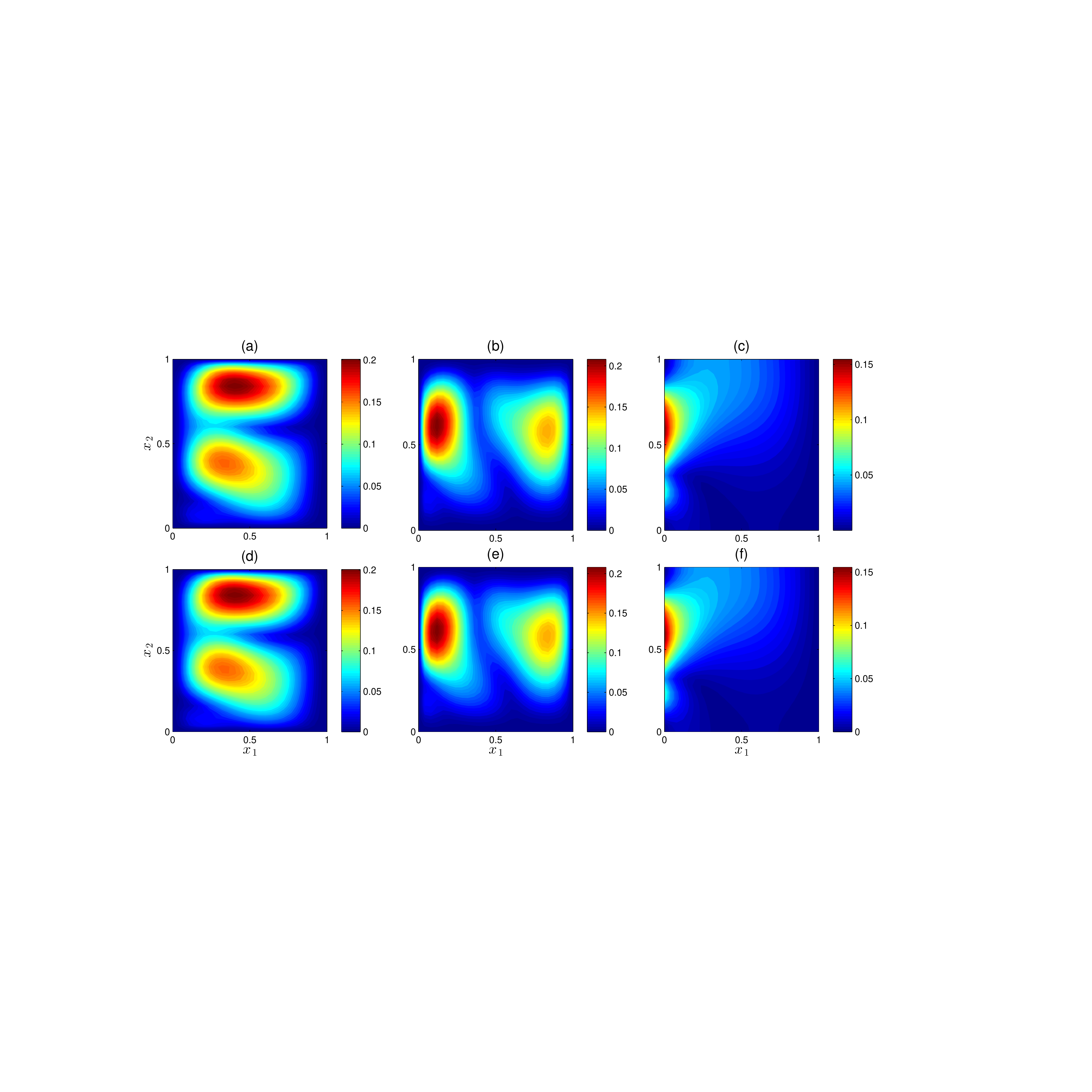}
    \caption{Standard deviation of solution fields obtained using both NISP method implementations. Subfigures (a, d), (b, e) and (c, f) correspond to $u_1$, $u_2$ and $T$ respectively, while (a, b, c) and (d, e, f) correspond to the reduced and standard NISP method implementations respectively.}
    \label{fig:fig8}
\end{figure}

\begin{table}[htbp]
\hypertarget{tab4}{}
\caption{Comparison of the approximation errors and computational costs
(seconds) obtained in the standard and reduced NISP method implementations for the Boussinesq flow problem. The
reduced dimensions and orders are also listed here.}
\begin{center}\scriptsize
\renewcommand{\arraystretch}{1.3}
\begin{tabular}{cc|ccc|cccccccc|c}
\toprule 
 &  & \multicolumn{3}{c}{Standard NISP} \vline & \multicolumn{8}{c}{Reduced NISP} \vline & \tabularnewline
$s_{1},s_{2}$ & {$p$} & {$Q$} & {$\varepsilon_{s}$} & {$\mathcal{C}_{s}$ } & {$d_{1}$} & {$\tilde{p}_{1}$} & {$\tilde{Q}_{1}$} & {$d_{2}$} & {$\tilde{p}_{2}$} & {$\tilde{Q}_{2}$} & {$\varepsilon_{r}$} & {$\mathcal{C}_{r}$ } & {$\mathcal{C}_{s} / \mathcal{C}_{r}$}\tabularnewline
\midrule
\midrule 
 & {$2$} & {$85$} & {$4.2\times10^{-4}$} & {$154$} & {$5$} & {$1$} & {$85$} & {$4$} & {$1$} & {$70$} & {$8.6\times10^{-3}$} & {$136$} & {$1.1$}\tabularnewline
\cmidrule{2-14}  
{$3$} & {$3$} & {$389$} & {$8.3\times10^{-5}$} & {$716$} & {$5$} & {$1$} & {$126$} & {$4$} & {$1$} & {$70$} & {$5.0\times10^{-3}$} & {$191$} & {$3.7$}\tabularnewline
\cmidrule{2-14}  
 & {$4$} & {$1457$} & {$2.9\times10^{-6}$} & {$2694$} & {$5$} & {$1$} & {$126$} & {$4$} & {$1$} & {$70$} & {$4.9\times10^{-3}$} & {$198$} & {$13.6$}\tabularnewline
\midrule 
 & {$2$} & {$145$} & {$4.7\times10^{-4}$} & {$260$} & {$6$} & {$1$} & {$145$} & {$5$} & {$1$} & {$126$} & {$1.2\times10^{-2}$} & {$242$} & {$1.1$}\tabularnewline
\cmidrule{2-14}  
{$4$} & {$3$} & {$849$} & {$3.9\times10^{-5}$} & {$1561$} & {$6$} & {$1$} & {$210$} & {$5$} & {$1$} & {$126$} & {$6.9\times10^{-3}$} & {$329$} & {$4.7$}\tabularnewline
\cmidrule{2-14}  
 & {$4$} & {$3937$} & {$3.0\times10^{-6}$} & {$7544$} & {$6$} & {$1$} & {$210$} & {$5$} & {$1$} & {$126$} & {$6.7\times10^{-3}$} & {$353$} & {$21.4$}\tabularnewline
\bottomrule
\end{tabular}
\par\end{center}
\end{table}

For various instances of $s_{1}$, $s_{2}$ and $p$, the approximation errors $\varepsilon_{s},\varepsilon_{r}$ and computational costs $\mathcal{C}_{s},\mathcal{C}_{r}$ (wall-times) observed in both
NISP algorithms, along with the reduced dimensions and orders observed at
the last iteration of the reduced NISP algorithm are listed in \hyperlink{tab4}{Table 4}. The number of iterations
was a constant and equal to $9$ for all cases, which once again indicates an invariance of 
the convergence rate with respect to $s_1$, $s_2$ and $p$.

Once again, due to the high degree of stochastic regularity exhibited
by the solutions, the approximation errors in the standard NISP method implementation
are observed to decay exponentially, while the predicted 
asymptotic upper bound for the approximation errors
in the reduced NISP method implementation is observed. The highest speedup 
factor observed is $\approx 20.4$.

\section{Summary and outlook}

A reduced NISP based method for uncertainty propagation in stochastic
multi-physics models is presented, demonstrated and compared against the standard NISP method using two numerical examples. At the expense of a
small approximation error, the reduced NISP method exhibits
signifcantly high speedup factors over the standard NISP method.
This is primarily due to
the much slower growth in the number of repeated runs of each
solver component in the former approach, which in turn can be attributed to the dimension
and order reduction steps for constructing an approximation of their respective input data.
Therefore, the curse of dimensionality,
which is the primary bottleneck for any spectral method in tackling
large multi-physics models with several independent stochastic components,
is demonstrably mitigated with the reduced NISP approach. 

In the intermediate dimension reduction step,
since the gPC coefficients of the local
variables are included in the matrix for which the SVD is computed, we
the local
stochastic dimension in each module would serve as a lower bound for the respective reduced dimension.
This suggests that the modular structure of the 
stochastic solver enables our proposed reduction strategy, which would 
fail to provide a reduced approximation if implemented
with a fully-coupled (monolithic) NISP solver. The lower bound on the reduced dimension can also be observed in the numerical experiments.

The two most significant overheads observed while
implementing the reduced NISP algorithm
were the computation of the Vandermonde matrix for the construction
of the reduced quadrature rule, and the computation of the global
gPC coefficient matrices from the reduced gPC coefficient matrices. Both these overheads can be
easily eliminated due to the embarrassingly parallel nature of the computations involved.
Moreover, in each numerical example, both the NISP algorithms were implemented using a sparse grid global quadrature
rule. If a full tensor grid were used, the computational
gains observed in the reduced NISP method implementation would be much higher. 
Furthermore, alternatives, for example, Schur complement based elimination, to the BGS partitioning approach can be explored
for reducing the overall computational costs.

The applicability of our proposed algorithm is limited to
models in which the global stochastic dimension is manageably low ($<20$), since the global gPC coefficients still need to be computed, stored and operated on. 
If the global stochastic dimension becomes prohibitively large
due to a particular module contributing a large number of independent uncertainties,
a Monte-Carlo based sampling approach can be employed in that module, while
other modules could still afford the use of spectral methods. Such a framework
has been recently demonstrated in [\hyperlink{ref41}{41}], and further explorations to develop a more general approach are currently
underway. Moreover,
for tackling models in which the solution regularity is low,
the proposed method can also be easily adapted towards multi-element gPC
[\hyperlink{ref42}{42}] and discontinuous wavelet [\hyperlink{ref43}{43}] based spectral representations. 

Extending the framework of dimension and order reduction to intrusive
spectral projection (ISP)- based uncertainty propagation would be particularly important in mitigating the overall intrusiveness and therefore, development
costs incurred by the solver modules. Moreover, derivative or active-subspace based dimension reduction
methods [\hyperlink{ref44}{44}] have provided potential alternative approaches which
are not yet fully explored. 
Futhermore, the additional complexities
and approximation errors that arise in tackling unsteady multi-physics models using spectral methods have led to several challenges and opportunities for improvement.

\section*{Acknowledgement}
This research was funded by the US Department of Energy, Office of
Advanced Computing Research and Applied Mathematics Program and partially
funded by the US Department of Energy NNSA ASC Program.

\newpage
\section*{Appendix A - Karhunen-Loeve expansion}
\hypertarget{appA}{}
\subsubsection*{Exponential covariance kernel}

Given an $n$-dimensional spatial domain $\Omega\subseteq\mathbb{R}^{n}$,
let $C_{u}:\Omega\times\Omega\rightarrow\mathbb{R}^{+}$ denote the exponential
covariance kernel of a spatially varying random field $u\in\Omega\rightarrow\mathbb{R}$.
Therefore, $C_{u}:\forall\boldsymbol{x}=\left[\begin{array}{ccc}
x_{1} & \cdots & x_{n}\end{array}\right]^{\mathbf{T}},\boldsymbol{y}=\left[\begin{array}{ccc}
y_{1} & \cdots & y_{n}\end{array}\right]^{\mathbf{T}}\in\Omega$, 
\begin{align*}
C_{u}\left(\boldsymbol{x},\boldsymbol{y}\right)=\exp\left(-\frac{\left\Vert \boldsymbol{x}-\boldsymbol{y}\right\Vert _{1}}{l}\right)=\prod_{j=1}^{n}\exp\left(-\frac{\left|x_{j}-y_{j}\right|}{l}\right),\tag{A.1}
\end{align*}
where $l$ denotes the correlation length. Subsequently, we can define
the KL expansion of $u$ using an infinite set of  random variables
$\left\{ \xi_{\boldsymbol{j}}:\boldsymbol{j}\in\mathbb{N}^{n}\right\}$, each having zero mean and unit variance, as follows. $\forall\boldsymbol{x}\in\Omega$, 

\hypertarget{eqA2}{}
\begin{align*}
u\left(\boldsymbol{x}\right)-\bar{u}\left(\boldsymbol{x}\right) & =\sum_{\boldsymbol{j}\in\mathbb{R}^{n}}\gamma_{\boldsymbol{j}}\left(\boldsymbol{x}\right)\xi_{\boldsymbol{j}}\\
 & =\sum_{j_{1}\in\mathbb{R}}\cdots\sum_{j_{n}\in\mathbb{R}}\gamma_{j_{1}\ldots j_{n}}\left(\boldsymbol{x}\right)\xi_{j_{1}\ldots j_{n}}\\
 & =\sum_{j_{1}\in\mathbb{R}}\cdots\sum_{j_{n}\in\mathbb{R}}\prod_{k=1}^{n}g_{j_{k}}\left(x_{k}\right)\xi_{j_{1}\ldots j_{n}}\tag{A.2}
\end{align*}
where $\forall j>0$, if $\zeta_{j}$ solves
\[
l\zeta_{j}+\tan\left(\frac{\zeta_{j}}{2}\right)=0,\tag{A.3}
\]
and $\zeta_{j+1}>\zeta_{j}>0$, then $\forall x\in\mathbb{R}$, 
\[
g_{j}\left(x\right)=\begin{cases}
2{\displaystyle \sqrt{\frac{l\zeta_{j}}{1+l^{2}\zeta_{j}^{2}}}\frac{\cos\left(\zeta_{j}x\right)}{\sqrt{\zeta_{j}+\sin\left(\zeta_{j}\right)}}} & j\ \mathrm{is\ odd},\\
2{\displaystyle \sqrt{\frac{l\zeta_{j}}{1+l^{2}\zeta_{j}^{2}}}\frac{\sin\left(\zeta_{j}x\right)}{\sqrt{\zeta_{j}-\sin\left(\zeta_{j}\right)}}} & j\ \mathrm{is\ even}.
\end{cases}\tag{A.4}
\]
Therefore, as is required in \hyperlink{sec41}{\S4.1} and \hyperlink{sec42}{\S4.2}, a truncated KL expansion
can be easily obtained from the single index form of the expansion
in \hyperlink{eqA2}{Eq. A.2}.

\renewcommand\thefigure{B\arabic{figure}}    
\setcounter{figure}{0}    
\newpage
\section*{Appendix B - Verification of modular solvers}
\subsection*{B1 - Poisson problem}
\hypertarget{appB1}{}
Given stochastic diffusion coefficients $a_{1}$ and $a_{2}$,
we choose analytical functions $u_{1}^{*}$, $u_{2}^{*}:\forall\boldsymbol{\xi}\in\Xi$,
\begin{align*}
u_{1}^{*}\left(\boldsymbol{x},\boldsymbol{\xi}\right) & =\frac{1}{\pi^{2}}\cos\left(\frac{\pi}{2}x_{1}\right)\sin\left(\pi x_{2}\right), &\boldsymbol{x}\in\Omega_{1},\\
u_{2}^{*}\left(\boldsymbol{x},\boldsymbol{\xi}\right) & =\frac{1}{\pi^{2}}\cos\left(\frac{\pi}{2}x_{1}\right)\sin\left(\pi x_{2}\right), &\boldsymbol{x}\in\Omega_{2},\tag{B.1}
\end{align*}
which solve the modified coupled-domain Poisson equations:$\forall\boldsymbol{\xi}\in\Xi$,
\begin{align*}
\boldsymbol{\nabla}^{\mathbf{T}}\left(a_{1}\left(\boldsymbol{x},\boldsymbol{\xi}_{1}\right)\boldsymbol{\nabla}u_{1}\left(\boldsymbol{x},\boldsymbol{\xi}\right)\right)+f_{1}^{*}\left(\boldsymbol{x},\boldsymbol{\xi}_{1}\right) & =0, & \boldsymbol{x}\in\Omega_{1},\\
\boldsymbol{\nabla}^{\mathbf{T}}\left(a_{2}\left(\boldsymbol{x},\boldsymbol{\xi}_{2}\right)\boldsymbol{\nabla}u_{2}\left(\boldsymbol{x},\boldsymbol{\xi}\right)\right)+f_{2}^{*}\left(\boldsymbol{x},\boldsymbol{\xi}_{2}\right) & =0, & \boldsymbol{x}\in\Omega_{2},\tag{B.2}
\end{align*}

Here, $\forall i\in\left\{ 1,2\right\} ,f_{i}^{*}:\forall\boldsymbol{x}\in\Omega_{i},\boldsymbol{\xi}_{i}\in\Xi_{i}$,
\begin{align*}
f_{i}^{*}\left(\boldsymbol{x},\boldsymbol{\xi}_{i}\right) & =\frac{5}{4}a_{i}\left(\boldsymbol{x},\boldsymbol{\xi}_{i}\right)\cos\left(\frac{\pi}{2}x_{1}\right)\sin\left(\pi x_{2}\right)+\frac{1}{2\pi}\frac{\partial a_{i}}{\partial x_{1}}\left(\boldsymbol{x},\boldsymbol{\xi}_{i}\right)\sin\left(\frac{\pi}{2}x_{1}\right)\sin\left(\pi x_{2}\right)\\
 & -\frac{1}{\pi}\frac{\partial a_{i}}{\partial x_{2}}\left(\boldsymbol{x},\boldsymbol{\xi}_{i}\right)\cos\left(\frac{\pi}{2}x_{1}\right)\cos\left(\pi x_{2}\right).\tag{B.3}
\end{align*}

\begin{figure}
    \hypertarget{figB1}
    \centering
    \includegraphics[bb=40bp 220bp 550bp 590bp,clip,scale=0.55]{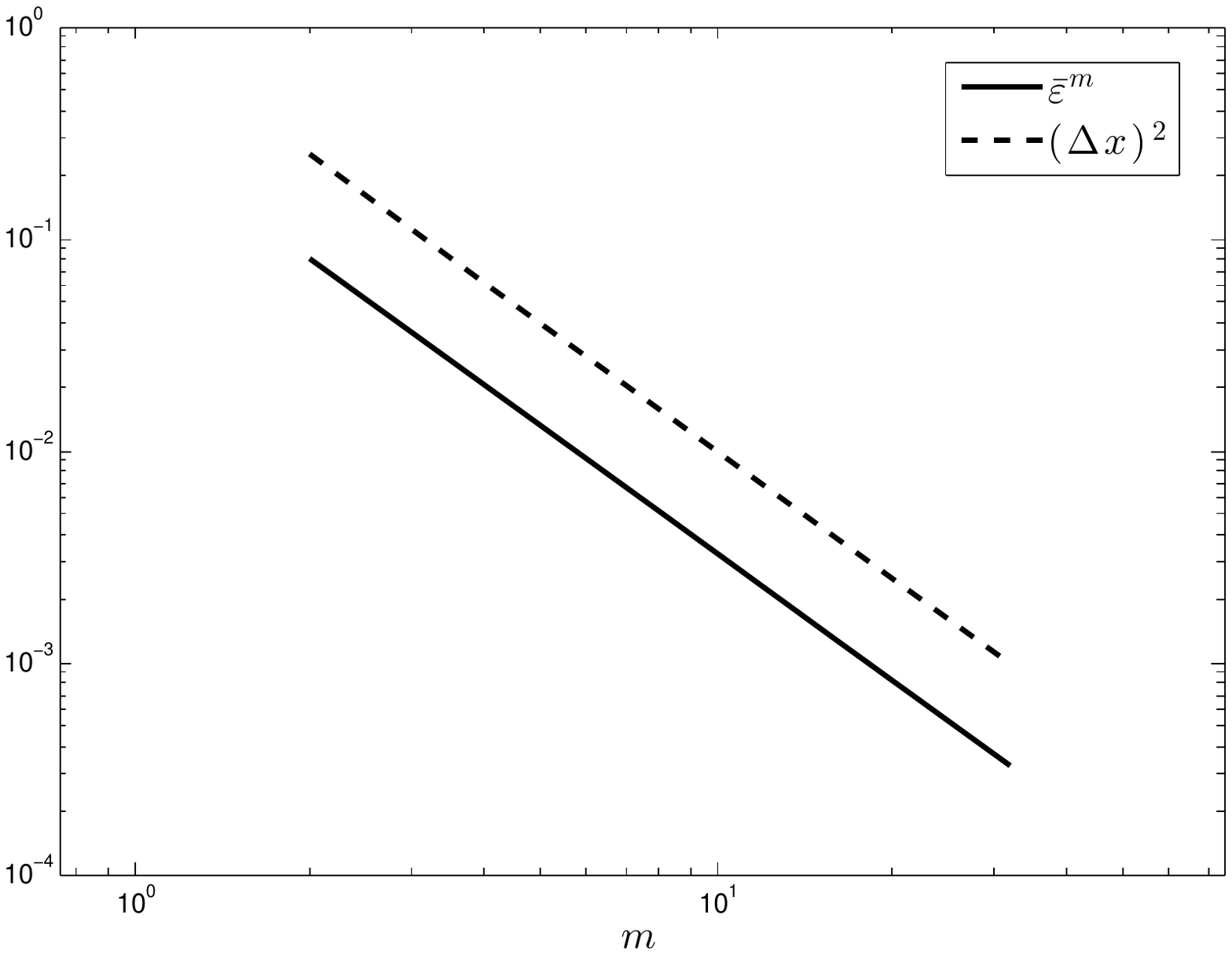} 
    \caption{Average error $\bar{\varepsilon}^{m}$ v.s.
$m$. $\Delta x$ denotes the node spacing.}
    \label{fig:figB1}
\end{figure}

If $u_{1}^{m},u_{2}^{m}$ denote the corresponding
finite element solutions obtained using $m\times m$ nodes, then $\varepsilon^{m}\left(\boldsymbol{\xi}\right):\forall\boldsymbol{\xi}\in\Xi$,
\begin{align*}
\varepsilon^{m}\left(\boldsymbol{\xi}\right) & =\sqrt{\frac{\sum_{i=1}^{2}\int_{\Omega_{i}}\left(u_{i}^{m}\left(\boldsymbol{x},\boldsymbol{\xi}\right)-u_{i}^{*}\left(\boldsymbol{x},\boldsymbol{\xi}\right)\right)^{2}d\boldsymbol{x}}{\sum_{i=1}^{2}\int_{\Omega_{i}}u_{i}^{*}\left(\boldsymbol{x},\boldsymbol{\xi}\right)^{2}d\boldsymbol{x}}}\tag{B.4}
\end{align*}
denotes the mean-square error at discretization level $m$. We used
a tolerance $\epsilon_{\mathrm{BGS}}=10^{-6}$ and set the the stochastic
dimensions to $s_{1}=s_{2}=4$. Subsequently, for different values
of $m$, we computed the sample average of $\varepsilon^{m}$,
denoted here as $\bar{\varepsilon}^{m}$, using $100$
Monte-Carlo solution samples. \hyperlink{figB1}{Figure B1} illustrates the expected second order
rate of decay in $\bar{\varepsilon}^{m}$.

\subsection*{B2 - Boussinesq flow problem}
\hypertarget{appB2}{}
Following the same procedure in \hyperlink{appB1}{Appendix B1}, we choose
analytical functions $u_{1}^{*},u_{2}^{*},p^{*},t^{*}:\forall\boldsymbol{x}\in\Omega,\boldsymbol{\xi}\in\Xi$,
\begin{align*}
u_{1}^{*}\left(\boldsymbol{x},\boldsymbol{\xi}\right) & =-\sin^{2}\left(\pi x_{1}\right)\sin\left(2\pi x_{2}\right),\\
u_{2}^{*}\left(\boldsymbol{x},\boldsymbol{\xi}\right) & =\sin\left(2\pi x_{1}\right)\sin^{2}\left(\pi x_{2}\right),\\
p^{*}\left(\boldsymbol{x},\boldsymbol{\xi}\right) & =\cos\left(\pi x_{1}\right)\cos\left(\pi x_{2}\right),\\
t^{*}\left(\boldsymbol{x},\boldsymbol{\xi}\right) & =\cos\left(\frac{\pi}{2}x_{1}\right)T_{h}\left(x_{2},\boldsymbol{\xi}_{2}\left(\boldsymbol{x},\boldsymbol{\xi}\right)\right),\tag{B.5}
\end{align*}
which solve the modified Boussinesq equations: $\forall\boldsymbol{x}\in\Omega,\boldsymbol{\xi}\in\Xi$,
\begin{align*}
\boldsymbol{\nabla}^{\mathbf{T}}\boldsymbol{u}\left(\boldsymbol{x},\boldsymbol{\xi}\right) & =0,\\
\left(\boldsymbol{u}\left(\boldsymbol{x},\boldsymbol{\xi}\right)^{\mathbf{T}}\boldsymbol{\nabla}\right)\boldsymbol{u}\left(\boldsymbol{x},\boldsymbol{\xi}\right)+\boldsymbol{\nabla}p\left(\boldsymbol{x},\boldsymbol{\xi}\right)\\
-\mathrm{Pr}\boldsymbol{\nabla}^{\mathbf{T}}\boldsymbol{\nabla}\boldsymbol{u}\left(\boldsymbol{x},\boldsymbol{\xi}\right)-\mathrm{Pr}\mathrm{Ra}\left(\boldsymbol{x},\boldsymbol{\xi}_{1}\left(\boldsymbol{x},\boldsymbol{\xi}\right)\right)T\left(\boldsymbol{x},\boldsymbol{\xi}\right)\boldsymbol{e}_{2}+\boldsymbol{f}_{u}^{*}\left(\boldsymbol{x},\boldsymbol{\xi}\right) & =\boldsymbol{0},\\
\left(\boldsymbol{u}\left(\boldsymbol{x},\boldsymbol{\xi}\right)^{\mathbf{T}}\boldsymbol{\nabla}\right)T\left(\boldsymbol{x},\boldsymbol{\xi}\right)-\boldsymbol{\nabla}^{\mathbf{T}}\boldsymbol{\nabla}T\left(\boldsymbol{x},\boldsymbol{\xi}\right)+f_{T}^{*}\left(\boldsymbol{x},\boldsymbol{\xi}\right) & =0.\tag{B.6}
\end{align*}

Here, $\forall\boldsymbol{x}\in\Omega,\boldsymbol{\xi}\in\Xi$, 
\begin{align*}
\boldsymbol{f}_{u}^{*}\left(\boldsymbol{x},\boldsymbol{\xi}\right) & =2\pi\left[\begin{array}{c}
\sin\left(2\pi x_{1}\right)\sin^{2}\left(\pi x_{1}\right)\sin^{2}\left(2\pi x_{2}\right)\\
\sin\left(2\pi x_{2}\right)\sin^{2}\left(\pi x_{2}\right)\sin^{2}\left(2\pi x_{1}\right)
\end{array}\right]\\
 & -\pi\left[\begin{array}{c}
\sin\left(2\pi x_{1}\right)\sin^{2}\left(\pi x_{1}\right)\sin^{2}\left(2\pi x_{2}\right)-\sin\left(\pi x_{1}\right)\cos\left(\pi x_{2}\right)\\
\sin\left(2\pi x_{2}\right)\sin^{2}\left(\pi x_{2}\right)\sin^{2}\left(2\pi x_{1}\right)-\sin\left(\pi x_{2}\right)\cos\left(\pi x_{1}\right)
\end{array}\right]\\
 & +2\pi^{2}\mathrm{Pr}\left[\begin{array}{c}
\cos\left(2\pi x_{1}\right)\sin\left(2\pi x_{2}\right)-2\sin^{2}\left(\pi x_{1}\right)\sin\left(2\pi x_{2}\right)\\
-\cos\left(2\pi x_{2}\right)\sin\left(2\pi x_{1}\right)+2\sin^{2}\left(\pi x_{2}\right)\sin\left(2\pi x_{1}\right)
\end{array}\right]\\
 & +\Pr\mathrm{Ra}\left(\boldsymbol{x},\boldsymbol{\xi}_{1}\left(\boldsymbol{\xi}\right)\right)\cos\left(\frac{\pi}{2}x_{1}\right)T_{h}\left(x_{2},\boldsymbol{\xi}_{2}\left(\boldsymbol{\xi}\right)\right)\boldsymbol{e}_{2},\\
f_{T}^{*}\left(\boldsymbol{x},\boldsymbol{\xi}\right) & =-\frac{\pi}{2}\left(\sin^{2}\left(\pi x_{1}\right)\sin\left(2\pi x_{2}\right)\sin\left(\frac{\pi}{2}x_{1}\right)+\frac{\pi}{2}\cos\left(\frac{\pi}{2}x_{1}\right)\right)T_{h}\left(x_{2},\boldsymbol{\xi}_{2}\left(\boldsymbol{\xi}\right)\right)\\
 & -\sin\left(2\pi x_{1}\right)\cos\left(\frac{\pi}{2}x_{1}\right)\sin^{2}\left(\pi x_{2}\right)\frac{\partial T_{h}}{\partial x_{2}}\left(x_{2},\boldsymbol{\xi}_{2}\left(\boldsymbol{\xi}\right)\right)\\
 & +\cos\left(\frac{\pi}{2}x_{1}\right)\frac{\partial^{2}T_{h}}{\partial x_{2}^{2}}\left(x_{2},\boldsymbol{\xi}_{2}\left(\boldsymbol{\xi}\right)\right).\tag{B.7}
\end{align*}

Let $\varepsilon^{m}\left(\boldsymbol{\xi}\right):\forall\boldsymbol{\xi}\in\Xi$,
\[
\varepsilon^{m}\left(\boldsymbol{\xi}\right)=\sqrt{\frac{\int_{\Omega}\left\Vert \left[\begin{array}{c}
u_{1}^{*}\left(\boldsymbol{x},\boldsymbol{\xi}\right)-u_{1}^{m}\left(\boldsymbol{x},\boldsymbol{\xi}\right) \\
u_{2}^{*}\left(\boldsymbol{x},\boldsymbol{\xi}\right)-u_{2}^{m}\left(\boldsymbol{x},\boldsymbol{\xi}\right) \\
p^{*}\left(\boldsymbol{x},\boldsymbol{\xi}\right)-p^{m}\left(\boldsymbol{x},\boldsymbol{\xi}\right) \\
t^{*}\left(\boldsymbol{x},\boldsymbol{\xi}\right)-t^{m}\left(\boldsymbol{x},\boldsymbol{\xi}\right)\end{array}\right]\right\Vert_{2}^{2}d\boldsymbol{x}}{\int_{\Omega}\left\Vert \left[\begin{array}{c}
u_{1}^{*}\left(\boldsymbol{x},\boldsymbol{\xi}\right) \\
u_{2}^{*}\left(\boldsymbol{x},\boldsymbol{\xi}\right) \\
p^{*}\left(\boldsymbol{x},\boldsymbol{\xi}\right) \\
t^{*}\left(\boldsymbol{x},\boldsymbol{\xi}\right)\end{array}\right]\right\Vert_{2}^{2}d\boldsymbol{x}}}
\]
denote the mean-square error between the exact and approximate solutions
$u_{1}^{m},u_{2}^{m},p^{m},t^{m}$,
computed using $m\times m$ cells.

Subsequently, the average error $\bar{\varepsilon}^{m}$,
was computed using 100 Monte-Carlo samples for various values of $m$, shown in \hyperlink{figB2}{Figure B2}, keeping the tolerance at $\epsilon_{\mathrm{BGS}}=10^{-6}$
and stochastic dimensions to $s_{1}=s_{2}=4$. As expected, a second order rate of decay
rate is observed in $\bar{\varepsilon}^{m}$. 

\begin{figure}
    \hypertarget{figB2}
    \centering
    \includegraphics[bb=40bp 220bp 550bp 590bp,clip,scale=0.55]{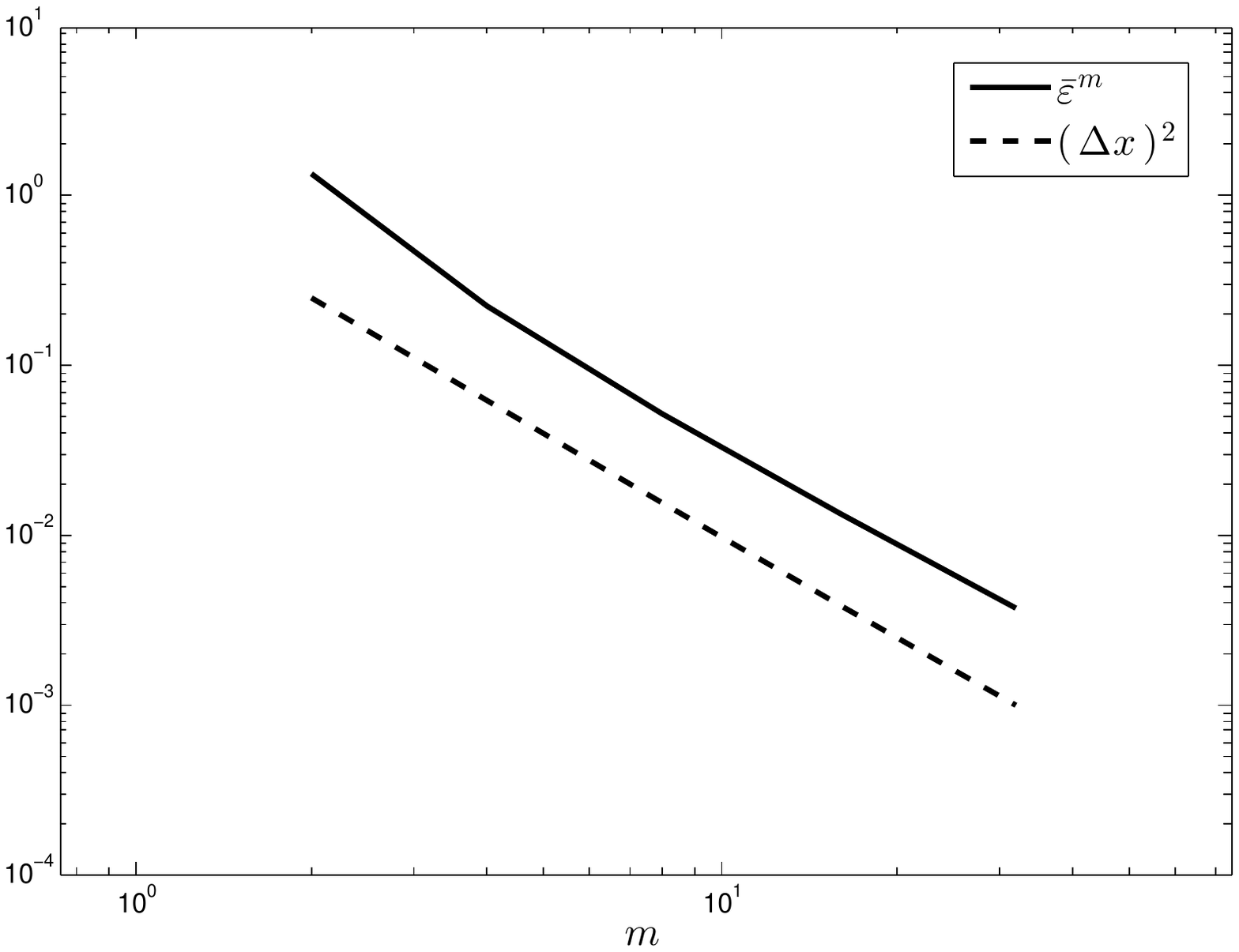} 
    \caption{Average error $\bar{\varepsilon}^{m}$ v.s.
$m$. $\Delta x$ denotes the cell side length.}
    \label{fig:figB2}
\end{figure}

\end{document}